\documentclass[unsortedaddress,nopreprint,twocolumn,numbers]{jasatex}

\usepackage{amssymb, latexsym}
\usepackage{amsmath} 
\usepackage{amsthm}
\usepackage{bm}
\usepackage{graphicx}
\usepackage{epstopdf}  
\usepackage[mathscr]{eucal}
\usepackage{enumerate}
 \graphicspath{{./Figures/}}

\newcommand{\be}{\begin{equation}}
\newcommand{\ee}{\end{equation}}
\renewcommand{\Vec}[1]{\boldsymbol{#1}}
\newcommand{\Mat}[1]{\boldsymbol{#1}}

\begin{document}

\title[Compressive beamforming]{Multiple and single  snapshot compressive beamforming}

\author{Peter Gerstoft}
\email[Corresponding author. Electronic mail: ]{gerstoft@ucsd.edu}
\affiliation{Scripps Institution of Oceanography, University of California San Diego, La Jolla, California 92093--0238}
\author{Angeliki Xenaki}
\affiliation{Department of Applied Mathematics and Computer Science, Technical University of Denmark, Kgs.Lyngby, 2800 Denmark}
\author{Christoph F. Mecklenbr\"auker}
\affiliation{Christian Doppler Lab, Inst.\ of Telecommunications, TU Wien, Gusshausstr.\ 25/389, 1040 Vienna, Austria}

\date{\today}

\begin{abstract}
For a  sound field observed on a sensor array, compressive sensing (CS) reconstructs the direction-of-arrival (DOA) of multiple sources using a sparsity constraint. The DOA estimation is posed as an underdetermined problem by expressing the acoustic pressure at each sensor as a phase-lagged superposition of source amplitudes at all hypothetical DOAs. Regularizing with an $\ell_1$-norm constraint renders the problem solvable with convex optimization, and  promoting sparsity gives high-resolution DOA maps. 
Here, the sparse source distribution is derived using maximum a posteriori (MAP) estimates for both  single and multiple snapshots. CS does not require  inversion of the data covariance matrix and thus works well even for a single snapshot where it gives  higher resolution than conventional beamforming. For multiple snapshots, CS outperforms conventional high-resolution methods, even with coherent arrivals and at low signal-to-noise ratio.
The superior resolution of CS is demonstrated with vertical array data from the SWellEx96 experiment for coherent multi-paths.
\end{abstract}


\maketitle

\section{\label{sec:intro}Introduction}

Direction-of-arrival (DOA) estimation refers to the localization of several sources from noisy measurements of the wavefield with an array of sensors.  DOA estimation can be expressed as a linear underdetermined problem with a sparsity constraint enforced on its  solution. The compressive sensing\cite{EladBook:2010, FoucartBook:2013} (CS) framework asserts that this is solved efficiently with a convex optimization procedure that promotes sparse solutions.

In DOA estimation, CS achieves high-resolution acoustic imaging\cite{MalioutovDOA:2005,XenakiCS:2014,Xenaki:2015}, outperforming traditional methods\cite{KrimDOA:1996}. Unlike the high-resolution subspace-based DOA estimators \cite{VanTreesBook, SchmidtMUSIC:1986}, DOA estimation via CS is reliable even with a single snapshot\cite{EdelmannCSDOA:2011,Mecklenbrauker:2013,FortunatiDOA:2014}. 

The  least absolute shrinkage and selection operator (LASSO) \cite{TibshiraniLasso1996} has been extended to multiple measurement vectors (here multiple snapshots) \cite{MalioutovDOA:2005,Wipf:2007}. 
They modify the LASSO objective function by introducing a mixed-norm penalty term 
that promotes spatial sparsity.
More specifically, the snapshots are combined with the $\ell_2$-norm, whereas the spatial samples are combined with the $\ell_1$-norm.
Multiple-snapshot CS offers several benefits over other high-resolution DOA estimators \cite{MalioutovDOA:2005,Wipf:2007,XenakiCS:2014}: 
1) It handles partially coherent arrivals.
2) It can be formulated with any number of snapshots, in contrast to, e.g., the Minimum Variance Distortion-free Response (MVDR) beamformer. 
3) Its flexibility in formulation enables extensions to sequential processing, and online algorithms \cite{Mecklenbrauker:2013}. 
Here, we show that CS achieves higher resolution than MUSIC and MVDR, even in scenarios that favor these classical high-resolution methods.  

In ocean acoustics, CS has found several applications in matched field processing\cite{MantzelCMFP:2012,Forero:2014}
 and in coherent passive fathometry for inferring sediment interfaces depths and their number\cite{YardimCSfathometer:2014}.
Various wave propagation phenomena from a single source (refraction, diffraction, scattering, ducting, reflection) lead to multiple partially coherent arrivals received by the array.
High-resolution beamformers cannot resolve these coherent arrivals.

CS for single snapshot has high-resolution capabilities and contrary to eigenvalue-based beamformers works for coherent arrivals\cite{MalioutovDOA:2005,XenakiCS:2014,FortunatiDOA:2014}.
CS is limited by  basis mismatch\cite{ChiBasisMismatch:2011} which occurs when the DOAs do not coincide with the look directions  of the angular spectrum, and by basis coherence.
Solutions to basis mismatch involve for example using atomic norm and solving the dual problem \cite{Candes:2013,Xenaki:2015} that are not addressed here.
%
Grid refinement alleviates basis mismatch for high signal to noise ratio (SNR) at the expense of increased computational complexity.
A denser grid causes increased coherence among the
steering vectors (basis coherence) which translates to bias and spread in the DOA estimates as demonstrated here. This is especially true in large two-dimensional or three-dimensional geo-acoustic inversion problems as e.g. seismic imaging\cite{YaoCSseismic:2011,YaoCSseismic:2013, Fan:2014}.

We use least squares optimization with an $\ell_{1}$-norm regularization term, also known as the LASSO \cite{TibshiraniLasso1996}, to formulate the DOA estimation problem for single and multiple snapshots. The LASSO formulation complies with statistical models as it provides a  maximum a posteriori (MAP) estimate, assuming a Gaussian data likelihood  and a Laplacian prior distribution for the source acoustic pressure\cite{YuanBayes2005,ParkBayes2008} for both single (Sec.\ \ref{sec:MAP}) and multiple snapshots\cite{Wipf:2007}  (Sec.\ \ref{sec:MSCS}). The LASSO is known to be a convex minimization problem and solved efficiently by interior point methods. 
In the LASSO formulation, Sec.\ \ref{sec:LASSO-path}, the reconstruction  accuracy depends on the choice of the regularization parameter that controls the balance between the data fit and the  sparsity of the  solution. We indicate that the regularization parameter can be found from the properties of the LASSO path\cite{TibshiraniLassoPath2011,Panahi:2012}, i.e., the evolution of the LASSO solution versus  the regularization parameter.

The main focus of the paper is  on performance evaluation for single and multiple snapshots using both simulated (Sec.\ \ref{sec:perf}) and real data (Sec. \ref{sec:swellex96}). Other excellent papers\cite{FortunatiDOA:2014} have already performed performance evaluation 
for single snapshot, consistent with our simulations. We are not aware of performance evaluation for multiple snapshots. 
 
In the following,
the $\ell_{p}$-norm of a vector $\mathbf{x} \in\mathbb{C}^{N}$ is defined as $\lVert\textbf{x}\rVert_{p} = \left(\sum_{n=1}^{N}\lvert x_{n}\rvert^{p}\right)^{1/p}$. By extension, the $\ell_{0}$-norm is defined as $\lVert\textbf{x}\rVert_{0}=\sum_{n=1}^{N} 1_{x_{n}\neq0}$ and the $\ell_{\infty}$-norm as $\lVert\textbf{x}\rVert_{\infty}=\max\limits_{1\le n\le N} \lvert x_{n}\rvert$.
For a matrix $\mathbf{F}\in\mathbb{C}^{M\times L}$ the Frobenius norm $\lVert \cdot \rVert_{\mathcal{F}}$  is defined as $\lVert \mathbf{F} \rVert_{\mathcal{F}}^{2} = \sum_{i=1}^{M}\sum_{j=1}^{L} \lvert f_{i,j}  \rvert^{2}$.

\section{Single snapshot DOA estimation}

We assume plane wave propagation and narrowband processing with a known sound speed. We consider the one-dimensional problem with a uniform linear array (ULA) of sensors with the source location characterized by the DOA of the associated plane wave, $\theta \in [-90^{\circ}, 90^{\circ}]$, with respect to the array axis.
%
 The propagation delay from the $i$th potential source to each of the $M$ array sensors is described by the steering (or replica) vector,
\begin{equation}
\mathbf{a}(\theta_{i})=\frac{1}{\sqrt{M}}\left[ 1,   e^{j \frac{2\pi d}{\lambda} 1\sin\theta_{i}}, \ldots , e^{j \frac{2\pi d}{\lambda} (M-1)\sin\theta_{i}}\right]^{T},
\label{eq:SteeringVector} 
\end{equation}
 where $\lambda$ is the wavelength and $d$ the sensor spacing.

Discretizing the half-space of interest, $\theta \in [-90^{\circ}, 90^{\circ}]$, into $N$ angular directions the DOA estimation problem can be expressed  as a  source reconstruction problem with the linear model,
\begin{equation}
\mathbf{y} = \mathbf{A}\mathbf{x} + \mathbf{n},
\label{eq:NoisyData}
\end{equation}
 where $\mathbf{y}\in \mathbb{C}^{M}$ is the complex-valued data vector from the measurements at the $M$ sensors, $\mathbf{x}\in \mathbb{C}^{N}$ is the unknown vector of the complex source amplitudes at all $N$ directions on the angular grid of interest and $\mathbf{n}\in\mathbb{C}^{M}$ is the additive noise vector. The sensing matrix,
\begin{equation}
\mathbf{A}=[\mathbf{a}(\theta_{1}), \cdots, \mathbf{a}(\theta_{N})],
\label{eq:SensingMatrix}
\end{equation}
maps the hypothetical sources $\mathbf{x}$ to the observations $\mathbf{y}$ and has as columns the steering vectors, Eq.~\eqref{eq:SteeringVector}, at all look directions.


In the following, the noise is generated as independent and identically distributed (iid) complex Gaussian. The array signal-to-noise ratio (SNR) is  defined as,
\begin{equation}
\mathrm{SNR}=10\log_{10}\frac{{\mathrm E}\left\lbrace\lVert\mathbf{A}\mathbf{x}\rVert_{2}^{2}\right\rbrace}{{\mathrm E}\left\lbrace\lVert\mathbf{n}\rVert_{2}^{2}\right\rbrace} \quad (\mathrm{dB}).
\label{eq:ArraySNR}
\end{equation}
%
\subsection{Sparse reconstruction with compressive sensing}

The problem of DOA estimation is to recover the set of non-zero components in the source vector $\mathbf{x}\in \mathbb{C}^{N}$, given the sensing matrix $\mathbf{A}_{M\times N}$ and an observation vector $\mathbf{y}\in \mathbb{C}^{M}$. Even though there are only a few sources $K<M$ generating the acoustic field, we are interested in a fine resolution on the angular grid to achieve precise localization such that $M\ll N$ and the problem in Eq.~\eqref{eq:NoisyData} is underdetermined. A way to solve this ill-posed inverse problem is  constraining the possible solutions with prior information.

Traditional methods solve the problem in Eq.~\eqref{eq:NoisyData} by seeking the solution with the minimum $\ell_{2}$-norm which provides the best data fit  ($\ell_{2}$-norm regularized least squares),
\begin{equation}
\widehat{\mathbf{x}}_{\ell_{2}} (\mu)=\underset{\mathbf{x}\in\mathbb{C}^{N}}{\arg\min} \; \lVert \mathbf{y} - \mathbf{A}\mathbf{x} \rVert_{2}^{2} + \mu \lVert \mathbf{x} \rVert_{2}^{2}.
\label{eq:CS_l2}
\end{equation}
The regularization parameter, $\mu \geq 0 $, controls the relative importance between the data fit and the $\ell_{2}$-norm of the solution. The minimization problem in Eq.~\eqref{eq:CS_l2} is convex with analytic solution, $\mathbf{\widehat{x}}_{\ell_{2}}(\mu) = \mathbf{A}^{H}\left(\mathbf{A}\mathbf{A}^{H}+ \mu \mathbf{I}_{M}\right)^{-1}\mathbf{y}$, where $\mathbf{I}_{M}$ is the $M \times M$ identity matrix. However, it aims to minimize the energy of the source $\bf x$ through the $\ell_{2}$-norm regularization term rather than its sparsity, hence the resulting solution is non-sparse. 

Conventional beamforming (CBF)\cite{VanTreesBook} is related to the $\ell_{2}$ solution for large $\mu$. From Eq.~\eqref{eq:CS_l2}:
\begin{equation}
\mathbf{\widehat{x}}_{\text{CBF}} = \lim\limits_{\mu \to \infty} (\mu\widehat{\mathbf{x}}_{\ell_{2}} (\mu)) = \mathbf{A}^{H}\mathbf{y}.
\label{eq:CBF}
\end{equation}
In principle, CBF combines the sensor outputs coherently to enhance the source signal at a specific look direction from the ubiquitous noise. CBF is robust to noise but suffers from low resolution and the presence of sidelobes.

Since $\mathbf{x}$ is sparse (there are only $K\ll N$ sources), it is appropriate to seek for the solution with the minimum $\ell_{0}$-norm, which counts the number of non-zero entries in the vector, to find a sparse solution. However, the $\ell_{0}$-norm minimization problem is a non-convex combinatorial problem which becomes computationally intractable even for moderate dimensions. The breakthrough of CS\cite{EladBook:2010, FoucartBook:2013} came with the proof that for sufficiently sparse signals, $K\ll N$, and sensing matrices with sufficiently incoherent columns the $\ell_{0}$-norm minimization problem is equivalent (at least in the noiseless case) to its  convex relaxation, the $\ell_{1}$-norm minimization problem\cite{CandesCStutorial:2008, BaraniukCSNotes:2007}. By replacing the $\ell_{0}$-norm with the convex $\ell_{1}$-norm, the problem can be solved efficiently with convex optimization even for large dimensions\cite{CVX, CVXTutorial, BoydBook}.

\begin{figure} 
\centering
	\includegraphics[width=8.6cm]{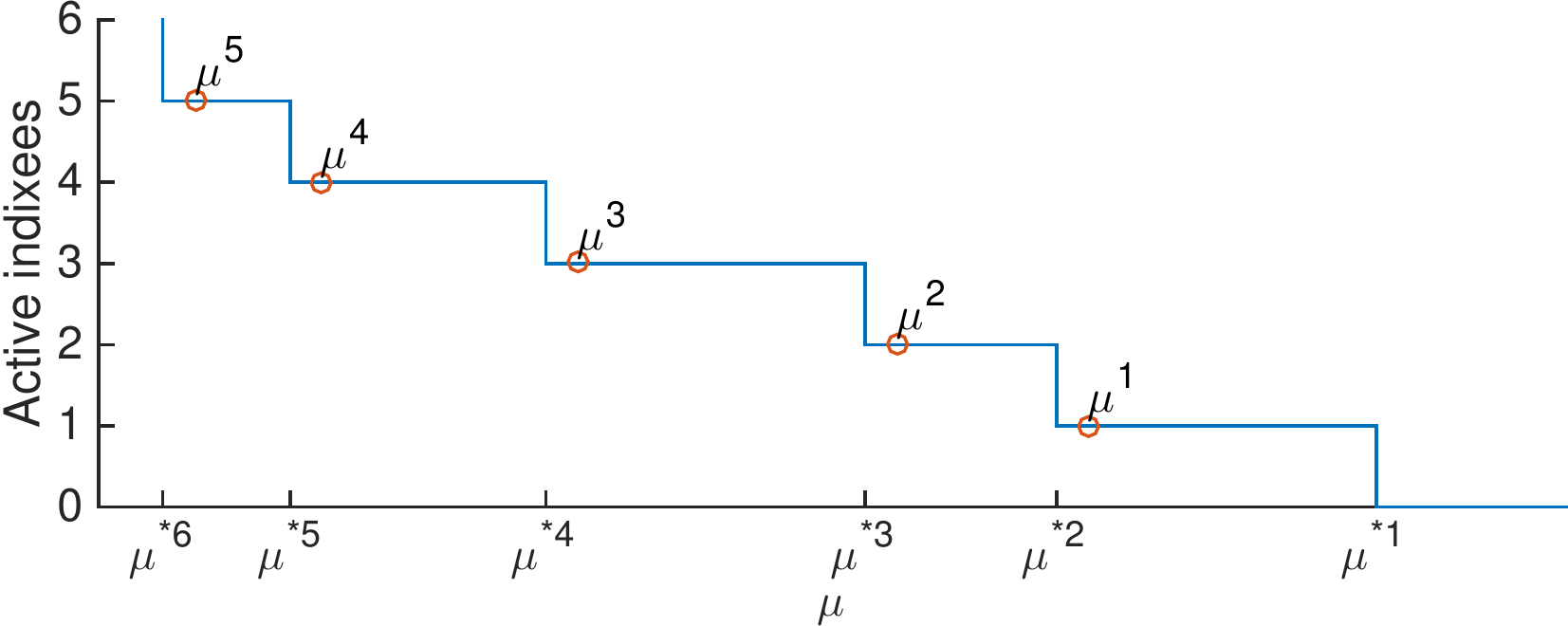}
\caption{Illustration of the LASSO path: Number of active indices versus the regularization parameter $\mu$. Increments in the active set occur at $\mu^{*p}$} \label{fig:lassopath}
\end{figure}

For noisy measurements, Eq.~\eqref{eq:NoisyData}, the $\ell_{1}$-norm minimization problem is formulated as 
\begin{equation}
\widehat{\mathbf{x}}_{\ell_{1}} (\epsilon)=\underset{\mathbf{x}\in\mathbb{C}^{N}}{\arg\min}\lVert\mathbf{x}\rVert_{1} \; \text{subject to} \;\lVert \mathbf{y} - \mathbf{A}\mathbf{x} \rVert_{2}\leq\epsilon,
\label{eq:CS_solution_noisy}
\end{equation}
where $\epsilon$ is the noise floor. The estimate $\widehat{\mathbf{x}}_{\ell_{1}} (\epsilon)$ has the minimum $\ell_{1}$-norm while it fits the data up to the noise level.
The problem in Eq.~\eqref{eq:CS_solution_noisy} can be equivalently written in an unconstrained form with the use of the regularizer $\mu\geq 0$, 
\begin{equation}
\widehat{\mathbf{x}}_{\ell_1}(\mu)=\underset{\mathbf{x}\in\mathbb{C}^{N}}{\arg\min}\; \lVert \mathbf{y} - \mathbf{A}\mathbf{x} \rVert_{2}^{2} + \mu \lVert \mathbf{x} \rVert_{1}.
\label{eq:CS_solution_lasso}
\end{equation}
The  sparse source reconstruction problem in Eq.~\eqref{eq:CS_solution_lasso} is a least squares optimization method regularized with the $\ell_{1}$-norm of the solution  $\mathbf{x}$ and provides the best data fit ($\ell_{2}$-norm term) for the sparsity level determined by  the regularization parameter $\mu$. The optimization problem in Eq.~\eqref{eq:CS_solution_lasso} is also known as the least absolute shrinkage and selection operator (LASSO) since the $\ell_{1}$ regularizer shrinks the coefficients of  $\mathbf{x}$ towards zero as the regularization parameter $\mu$ increases\cite{TibshiraniLasso1996}. This is  illustrated in Fig.\   \ref{fig:lassopath}. For every $\epsilon$ there exists a $\mu$ so that the estimates in Eq.~\eqref{eq:CS_solution_noisy} and Eq.~\eqref{eq:CS_solution_lasso} are equal.

Once the active DOAs are recovered, by solving Eq.~\eqref{eq:CS_solution_noisy} or equivalently Eq.~\eqref{eq:CS_solution_lasso}, the unbiased complex source amplitudes are determined from,
\begin{equation}
\widehat{\mathbf{x}}_{\text{CS}} = \mathbf{A}_{\mathcal{M}}^{+}\mathbf{y},
\label{eq:AmplitudeEstimation}
\end{equation}
where  $\mathbf{A}_{\mathcal{M}}\in\mathbb{C}^{N\times K}$ contains only the ``active'' steering vectors  associated with non-zero components in the solution $\widehat{\mathbf{x}}_{\ell_1}(\mu)$ and $\mathbf{A}_{\mathcal{M}}^{+}$ is its Moore--Penrose pseudoinverse.

For a given sparsity level $K$ and corresponding set of active indexes $\mathcal{M}$, i.e. $|\mathcal{M}|=K$, Eq.\ \eqref{eq:AmplitudeEstimation} finds the best data fit. Thus, if the active sensing matrix  $\mathbf{A}_{\mathcal{M}}$  has   sufficiently incoherent columns 
it represents the solution to the $\ell_0$ problem
\begin{equation}
\widehat{\mathbf{x}}_{\ell_{0}} (K)=\underset{\mathbf{x}\in\mathbb{C}^{N}}{\arg\min}\;\lVert \mathbf{y} - \mathbf{A}\mathbf{x} \rVert_{2}\ \text{subject to} \; \lVert\mathbf{x}\rVert_{0} =K
\label{eq:CS_ell0}
\end{equation}

\subsection{MAP estimate via LASSO}
\label{sec:MAP}

We use the LASSO formulation, Eq.~\eqref{eq:CS_solution_lasso}, to solve the DOA estimation problem in favor of sparse solutions. The choice of the (unconstrained) LASSO formulation over the constrained formulation, Eq.~\eqref{eq:CS_solution_noisy}, allows the sparse reconstruction method to be interpreted in a statistical Bayesian setting, where the unknowns $\mathbf{x}$ and the observations $\mathbf{y}$ are both treated as stochastic (random) processes, by imposing a prior distribution on the solution  $\mathbf{x}$ which promotes sparsity\cite{TibshiraniLasso1996,YuanBayes2005,ParkBayes2008}. 

 Bayes theorem connects the posterior distribution $p(\mathbf{x} | \mathbf{y} )$ of the model parameters $\mathbf{x}$ conditioned on the data $\mathbf{y}$, with the data likelihood $p(\mathbf{y}| \mathbf{x})$, the prior distribution of the model parameters $p(\mathbf{x})$ and the marginal distribution of the data $p(\mathbf{y})$,
\begin{equation}
p(\mathbf{x} | \mathbf{y} )= \frac{p(\mathbf{y}| \mathbf{x})p(\mathbf{x})}{p(\mathbf{y})}.
\label{eq:BayesRule}
\end{equation} 
Then, the maximum a posteriori (MAP) estimate is,
\begin{equation}
\begin{aligned}
\widehat{\mathbf{x}}_{\text{MAP}} &= \underset{\mathbf{x}}{\arg\max}\; \ln p(\mathbf{x} | \mathbf{y} )\\
&= \underset{\mathbf{x}}{\arg\max}\; \left[ \ln p(\mathbf{y}| \mathbf{x}) + \ln p(\mathbf{x})\right]\\
&= \underset{\mathbf{x}}{\arg\min}\; \left[ -\ln p(\mathbf{y}| \mathbf{x}) - \ln p(\mathbf{x})\right],
\end{aligned}
\label{eq:MAPestimate}
\end{equation}
where the marginal distribution of the data $p(\mathbf{y})$ is omitted since it is independent of the model $\mathbf{x}$.

Based on a complex Gaussian noise model with independent and identically distributed (iid) real and imaginary parts, $\mathbf{n} \sim {\cal CN}({\bf 0}, \sigma^2 {\bf I})$, the likelihood of the data is also complex Gaussian distributed $p(\mathbf{y}| \mathbf{x}) \sim {\cal CN}(\mathbf{A}\mathbf{x}, \sigma^2 {\bf I})$,
\begin{equation}
p(\mathbf{y}| \mathbf{x}) =\pi^{-N} \sigma^{-2N} {\rm e}^{ -\frac{\lVert{\mathbf y}-{\mathbf A}{\mathbf x}\rVert^{2}_{2}}{\sigma^{2}}}.
\label{eq:Likelihood}
\end{equation}
Following \cite{He-Xie-Ding-Cichocki2007}, we assume that the coefficients of the solution  $\mathbf{x}$ are iid and follow a Laplacian-like distribution (for complex random variables). Such a prior has been shown to encourage sparsity in many situations because of the heavy tails and sharp peak at zero.The corresponding prior is
\renewcommand{\Re}{\mathop{\mathrm{Re}}}
\renewcommand{\Im}{\mathop{\mathrm{Im}}}
\begin{equation}
p(\mathbf{x})\propto \prod\limits_{i=1}^{N} {\rm e}^{-\frac{\sqrt{(\Re x_i)^2+(\Im x_i)^2}}{\nu}} = {\rm e}^{-\frac{\lVert \mathbf{x}\rVert_{1}}{\nu}}.
\label{eq:LaplacePrior}
\end{equation}
The LASSO estimate, Eq.~\eqref{eq:CS_solution_lasso}, can be interpreted as the MAP estimate,
\begin{equation}
\widehat{\mathbf{x}}_{\text{MAP}} \!\!= \!\!\underset{\mathbf{x}}{\arg\min} \left[ \lVert{\mathbf y}-{\mathbf A}{\mathbf x}\rVert^{2}_{2} + \mu \lVert \mathbf{x}\rVert_{1}\right] \!\!=\!\!\widehat{\mathbf{x}}_{\ell_1}(\mu),
\label{eq:MAPlasso}
\end{equation}
where $\mu=\sigma^2/\nu$. 

Equation \eqref{eq:LaplacePrior} imposes no restriction on the source phases. Here, the phase is  assumed uniformly $[0, 2\pi )$ distributed.

\section{Multiple-snapshot DOA estimation} 
\label{sec:MSCS} 

Even though for moving sources it befits to solve one optimization problem for each snapshot sequentially\cite{Mecklenbrauker:2013}, for stationary scenarios, the sensor data statistics can be aggregated across snapshots to provide a more stable estimate. 
Multiple snapshots are  referred to as multiple measurement vectors and the recovery might have better performance than single measurement vectors\cite{Cotter:2005}. Potentially the recovery can be made more robust by using a likelihood function Eq. \eqref{eq:Likelihood} with colored noise (full covariance matrix) or based on the Huber norm\cite{Ollila2015}.
For the  multiple-snapshot case, all snapshot are collected into one matrix, 
\begin{equation}
\mathbf{Y} = \mathbf{A}\mathbf{X} + \mathbf{N},
\label{eq:MultiSnapNoisyData}
\end{equation}
where, for $L$ snapshots, $\mathbf{Y} = [\mathbf{y}(1),\cdots,\mathbf{y}(L)]$ and $\mathbf{N}= [\mathbf{n}(1),\cdots,\mathbf{n}(L)]$ are $M\times L$ matrices with the measurement and noise vectors per snapshot as columns, respectively, and $\mathbf{X}$ is the $N\times L$ signal with the complex source amplitudes at the $N$ look directions per snapshot as columns. For stationary sources the matrix $\mathbf{X}= [\mathbf{x}(1),\cdots,\mathbf{x}(L)]$ exhibits row sparsity, i.e., it has a constant sparsity profile for every column, since the few existing sources are associated with the same DOA for all snapshots.
As the sources are stationary it makes sense to sum the source energy across all snapshots, giving the row norm ${\bf x}^{\ell_2}$
\begin{equation}
{\bf x}^{\ell_2}= \left(\sum_{l=1}^L |{\bf X}_{ \centerdot l} |^2 \right)^{1/2} ~.
\end{equation}
This quantity is sparse and in analogy with the single snapshot case we impose  a Laplacian-like prior
%
\begin{equation}
p(\mathbf{X})=p(\mathbf{x}^{\ell_{2}})\propto\exp(-\lVert\mathbf{x}^{\ell_{2}}\rVert_{1}/\nu),   
\end{equation}
 with no phase assumption. Similar to Eq.\ \eqref{eq:LaplacePrior} we assume the phase is uniformly iid  distributed on $[0, 2\pi )$.
 
We assume an iid complex Gaussian distribution for the data likelihood 
\begin{equation}
p(\mathbf{Y} | \mathbf{X})\propto \exp(-\lVert \mathbf{Y} -\mathbf{A}\mathbf{X} \rVert_{\mathcal{F}}^2/\sigma^2).
\end{equation}
Using Bayes theorem, the MAP solution is then
\begin{equation}
\begin{aligned}
\widehat{\mathbf{X}}&={\arg\max}\; p(\mathbf{Y} | \mathbf{X})p(\mathbf{X})\\
&=\underset{\mathbf{X}\in\mathbb{C}^{N\times L}}{\arg\min}
 \; \lVert \mathbf{Y} - \mathbf{A}\mathbf{X} \rVert_{\mathcal{F}}^{2} + \mu \lVert \mathbf{x}^{\ell_{2}} \rVert_{1}.
\end{aligned}
\label{eq:CS_lasso_MultiSnap}
\end{equation}
In this formulation we search for a sparse solution via the $\ell_1$ constraint. The source amplitude can, however, vary across snapshot. This is in contrast to covariance-matrix based beamforming that just inverts for the average source power.
 The processing performance can be improved by doing an eigenvalue decomposition of $\mathbf{X}$ and retaining just the largest  eigenvalues; see Refs.\onlinecite{MalioutovDOA:2005,XenakiCS:2014}. The smaller eigenvalues contain mostly noise so this improves processing. However, this  eigenvalue decomposition is not done here as this has features similar to forming a sample covariance matrix.

Once the active steering vectors have been recovered, the unbiased source amplitudes are estimated for each snapshot, similar to the single snapshot case, Eq.\ \eqref{eq:AmplitudeEstimation},
\begin{equation}
\widehat{\mathbf{X}}_{\text{CS}} = \mathbf{A}_{\alpha}^{+}\mathbf{Y},
\label{eq:AmplitudeEstimationMMV}
\end{equation}
If desired, an average power estimate  $\mathbf{x}^{\ell_{2}}_{\text{CS}} $ can be  obtained from the $\ell_{2}$-norm of the rows of  $\widehat{\mathbf{X}}_{\text{CS}} $, with the $i$th element squared of  $\mathbf{x}^{\ell_{2}}_{\text{CS}}$ being the source power estimate at $\theta_i$.

For reference, the  CBF, MVDR, and MUSIC use the data sample covariance matrix,
\begin{equation}
\mathbf{C} = \frac{1}{L}\mathbf{Y}\mathbf{Y}^H.
\label{eq:SCM}
\end{equation}
The beamformer power for CBF and MVDR respectively is then,
\begin{alignat}{2}
P_{\rm CBF}(\theta)   &=\mathbf{w}_{\rm CBF}^H(\theta)    \mathbf{C}\mathbf{w}_{\rm CBF}(\theta) \\
P_{\rm MVDR}(\theta)&= \mathbf{w}_{\rm MVDR}^H(\theta)\mathbf{C}\mathbf{w}_{\rm MVDR}(\theta),
\label{eq:beam}
\end{alignat}
where the corresponding weight vectors are given by,
\begin{alignat}{2}
 \mathbf{w}_{\rm CBF}(\theta)&= \mathbf{a}(\theta) \\
 \mathbf{w}_{\rm MVDR}(\theta)&= \frac{\mathbf{C}^{-1}\mathbf{a}(\theta)}{\mathbf{a}^H(\theta)\mathbf{C}^{-1}\mathbf{a}(\theta)}.
 \label{eq:weigth}
\end{alignat}

The CBF can also be based directly on snapshots, as the single snapshot CBF Eq.~\eqref{eq:CBF} can be generalized to multiple snapshots,  $\mathbf{\widehat{X}}_{\text{CBF}} =\mathbf{A}^{H}\mathbf{Y}$.
The power estimates $P_{\rm CBF}(\theta_i)$, $P_{\rm MVDR}(\theta_i)$, and the corresponding $i$th squared component of $\mathbf{x}^{\ell_{2}}_{\text{CS}} $ are thus comparable. Note that since the MVDR weights in Eq.~\eqref{eq:weigth} involve the inverse of the sample covariance matrix, MVDR requires a full rank $\mathbf C$, i.e., $L\geq M$ snapshots.

%


The MUSIC\cite{VanTreesBook} is based on the eigendecomposition of the data sample covariance matrix Eq. \eqref{eq:SCM} and the separation of the signal and the noise subspace, 
\begin{equation}
\mathbf{C} = {\mathbf{U}}_{s}{\mathbf{\Lambda}}_{s}{\mathbf{U}}_{s}^{H} +\mathbf{U}_{n}\mathbf{\Lambda}_{n}\mathbf{U}_{n}^{H}.
\label{eq:Rsvd}
\end{equation}
The signal eigenvectors $\mathbf{U}_{s}$ corresponding to the largest eigenvalues, $\mathbf{\Lambda}_{s}$, are at the same subspace as the steering vectors, Eq.~\eqref{eq:SteeringVector}, while the noise eigenvectors $\mathbf{U}_{n}$ are orthogonal to the subspace of the steering vectors thus $\mathbf{a}(\theta)^{H}\mathbf{U}_{n} = \mathbf{0}$. 
MUSIC uses the orthogonality between the signal and the noise subspace to locate the maxima in the spectrum, 
\begin{equation}
P_{\rm MUSIC}(\theta) = \frac{1}{\mathbf{a}(\theta)^{H}\mathbf{U}_{n}\mathbf{U}_{n}^{H}\mathbf{a}(\theta)}.
\label{eq:Pmusic}
\end{equation}

Both MVDR and MUSIC overcome the resolution limit of the conventional beamformer by exploiting signal information conveyed by the data sample matrix. However, their performance depends on the eigenvalues of the data sample matrix thus it degrades with few snapshots, when the data sample matrix is rank deficient, and in the presence of coherent sources, when the signal subspace is reduced (Ch.~9 in Ref.\onlinecite{VanTreesBook}). CS does not have these limitations as it utilizes directly the measured pressure $\mathbf{Y}$.

%
%

\section{Regularization parameter selection}

The choice of the regularization parameter $\mu$ in Eq.~\eqref{eq:CS_solution_lasso}, also called the LASSO shrinkage parameter, is crucial as it controls the balance between the sparsity of the estimated solution and the data fit determining the quality of the reconstruction.

For large $\mu$, the solution is very sparse (with small $\ell_{1}$-norm) but the data fit is poor as indicated in Fig.\ \ref{fig:lassopath}. As $\mu$ decreases towards zero, the data fit is gradually improved since the corresponding solutions become less sparse. Note that for $\mu = 0$ the solution in Eq.~\eqref{eq:CS_solution_lasso} becomes the unconstrained least squares solution. Since the LASSO path is derived and demonstrated for a single observation, the statistics of the source signal or noise is irrelevant. 

\subsection{The LASSO path}
\label{sec:LASSO-path}

As the regularization parameter $\mu$ evolves from $\infty$ to $0$, the LASSO solution in Eq.~\eqref{eq:CS_solution_lasso} changes continuously following a piecewise smooth trajectory referred to as the solution path or the LASSO path\cite{TibshiraniLassoPath2011, Panahi:2012}.
 In this section, we show that the singularity points in the LASSO path are associated with a change in the sparsity of the solution and can be used to indicate an adequate choice for $\mu$.

We obtain the full solution path using convex optimization to solve Eq.~\eqref{eq:CS_solution_lasso} iteratively for different values of $\mu$. We use the CVX toolbox for disciplined convex optimization that is available in the Matlab environment. It uses interior point solvers to obtain the global solution of a well-defined optimization problem\cite{CVX, CVXTutorial, BoydBook}.

Let $L(\mathbf{x}, \mu)$ denote the objective function in Eq.~\eqref{eq:CS_solution_lasso},
\begin{equation}
L(\mathbf{x}, \mu) = \lVert \mathbf{y} - \mathbf{A}\mathbf{x} \rVert_{2}^{2} + \mu \lVert \mathbf{x} \rVert_{1}.
\label{eq:Lagrangian}
\end{equation}
\noindent The value $\widehat{\mathbf{x}}$ minimizing Eq.~\eqref{eq:Lagrangian} is found from its subderivative,
\begin{equation}
\begin{aligned}
&\partial_{\mathbf{x}} L(\mathbf{x}, \mu) = 2\mathbf{A}^{H}\left(\mathbf{A}\mathbf{x}-\mathbf{y}\right)  + \mu\partial_{\mathbf{x}} \lVert \mathbf{x} \rVert_{1},
\end{aligned}
\label{eq:DualFunction}
\end{equation}
where the subdifferential operator $\partial_{\mathbf{x}}$ is a generalization of the partial differential operator for functions that are not differentiable everywhere (Ref.\onlinecite{BoydBook} p.338). The subgradient for the $\ell_{1}$-norm is the set of vectors,
\begin{equation}
\partial_{\mathbf{x}} \lVert \mathbf{x} \rVert_{1} = \left\lbrace \mathbf{s} :\; \lVert \mathbf{s} \rVert_{\infty}\leq 1, \; \mathbf{s}^{H}\mathbf{x} = \lVert \mathbf{x}\rVert_{1} \right\rbrace,
\label{eq:L1normSubgradient}
\end{equation}
 which implies,
\begin{equation}
\begin{array}{ll}
s_{i} = \frac{x_{i}}{\lvert x_{i}\rvert}, &\;  x_{i}\neq 0 \\
\lvert s_{i} \rvert \le 1, &\; x_{i}= 0,
\end{array} 
\label{eq:ComponentwiseSubgradient}
\end{equation}
i.e., for every active element $x_{i}\neq 0$ of the  vector $\mathbf{x}\in\mathbb{C}^{N}$, the corresponding element of the subgradient is a unit vector in the direction of $x_{i}$. For every null element $x_{i}= 0$ the corresponding element of the subgradient has magnitude less than or equal to one. Thus, the magnitude of the subgradient is uniformly bounded by unity, $\lVert \mathbf{s} \rVert_{\infty}\leq 1$.

Denote,
\begin{equation}
\mathbf{r} = 2\mathbf{A}^{H}\left(\mathbf{y}- \mathbf{A}\widehat{\mathbf{x}} \right),
\label{eq:ResidualVector}
\end{equation}
the beamformed residual vector for the estimated solution $\widehat{\mathbf{x}}$. Since Eq.~\eqref{eq:Lagrangian} is convex, the global minimum is attained if $\mathbf{0} \in \partial_{\mathbf{x}} L(\mathbf{x}, \mu)$ which leads to the necessary and sufficient condition
\begin{equation}
 \mu^{-1}\mathbf{r} \quad\in\quad \partial_{\mathbf{x}} \lVert \mathbf{x} \rVert_{1}.
\label{eq:MinDualFunction}
\end{equation}
Then, from Eq.~\eqref{eq:ComponentwiseSubgradient} and Eq.~\eqref{eq:MinDualFunction}, the coefficients $r_{i}  =  2\mathbf{a}_{i}^{H}\left(\mathbf{y}- \mathbf{A}\widehat{\mathbf{x}} \right)$ of the beamformed residual vector $\mathbf{r} \in \mathbb{C}^{N}$ have amplitude such that,
\begin{equation}
\begin{array}{ll}
\lvert r_{i} \rvert =\mu, & \;  \widehat{x}_{i}\neq 0 \\
\lvert r_{i} \rvert \le \mu, & \; \widehat{x}_{i}= 0,
\end{array}
\label{eq:BoundedResidual}
\end{equation}
 i.e., whenever a component of $\widehat{\mathbf{x}}$ becomes non-zero, the corresponding element of the beamformed residual hits the boundary identified with the regularization parameter, $\lVert \mathbf{r} \rVert_{\infty}\leq \mu$. 

 For multiple snapshots, with the $ \widehat{\mathbf{X}}$ determined from Eq.\ \eqref{eq:CS_lasso_MultiSnap}, the beamformed residuals become
\be
\mathbf{R} = 2\mathbf{A}^{H}\left(\mathbf{Y}- \mathbf{A}\widehat{\mathbf{X}} \right), \quad \mathbf{r}_i =\sqrt{ \sum_{j=1}^L|\mathbf{R} _{ij} |^2  }~.
\label{eq:ResidualVectorMMV}
\ee 
The values of $\mu $ when changes in sparsity appear are obtained similarly to the single snapshot case.

\subsection{\label{sec:LassoAlgorithm}Algorithm for the LASSO path}

\begin{table}[h!]
\begin{center}
\caption{Fast iterative algorithm to solve the LASSO problem (\ref{eq:CS_solution_lasso}) for a desired sparsity level $K$ and estimating the unbiased complex source amplitudes  (\ref{eq:AmplitudeEstimation}).} 
\label{algo:fast}
\begin{tabular}{ll}
\hline\hline
    & Given: $\Mat{A}\in\mathbb{C}^{N\times M}$, $\Vec{y}\in\mathbb{C}^N$, $K\in{\mathbb N}$ ,  $F\in]0,1[$ \\ \hline 
1: &  Initialize   $i = 0$, $\Vec{x}_{\ell_1}^0 = \Vec{0} $, $\Vec{r}^{0}=2\Mat{A}^H  \Vec{y}  $  \\ \\  
2: &  \text{while} $|\mathcal{M}_{i}| < K$ \\
    & \hspace{2ex} $i = i + 1$ \\
3: & \hspace{2ex} $\mu^{i} =  (1-F)\mathop{\mathrm{peak}}\! \left(\Vec{r}^{i-1} , K \right)+ F\mathop{\mathrm{peak}}\! \left(\Vec{r}^{i-1} , K+1 \right) $\\
4a: &  \hspace{2ex} $\Vec{x}_{\ell_1}^{i} = \mathop{\mathrm{solution\ to\ Eq.\ (\ref{eq:CS_solution_lasso})\ for\ }} \Mat{A}, \Vec{y}, \mu = \mu^{i} $ \\
4b &  \hspace{2ex} $\Vec{r}^{i}=2\Mat{A}^H \left( \Vec{y} - \Mat{A}\Vec{x}_{\ell_1}^{i} \right)  $ \\
5: &  \hspace{2ex}    $\mathcal{M}_i = \{ m \Big|\, |x_{{\ell_1},m}^i| > \delta_i \}$,  $\delta_i= \epsilon \|\Vec{x}_{\ell_1}^i\|_{\infty} $  \\
  & \text{end} \\ \\
6:& \text{if} $|\mathcal{M}_{i}| > K$ \\
7:&  \hspace{2ex}    $\mathcal{M}_i = \{ m \Big|\, |x_{{\ell_1},m}^i| > \delta_i \}$,  $\delta_i= {\rm peak} (| \Vec{x}_{\ell_1}^i | ,K)$  \\
  & \text{end} \\
    & $\mathcal{M}  =\mathcal{M}_{i}  $ \\ \\
8: & $\widehat{\mathbf{x}}_{\ell_1}(\mu^i) = \Vec{x}_{\ell_1}^{i}$ \\
9: & $\widehat{\Vec{x}}_{\rm CS} =\Mat{A}_{\mathcal{M}}^+ \Vec{y} $ \\
10: & Output: $\mu^i$, $\widehat{\mathbf{x}}_{\ell_1}(\mu^i)$, $\widehat{\Vec{x}}_{\rm CS}$, $\mathcal{M}$. \\
\hline\hline
\end{tabular}
\end{center}
\end{table}

Although many algorithms exist for solving the LASSO problem, we have good experience with the  algorithm in Table \ref{algo:fast} as it is reasonable fast and accurate.
 Sec \ref{sec:LASSO-path} is used for formulating an algorithm where the values of $\mu$ for different sparsity levels are indicated by the dual solution $\mathbf r$, solving the dual problem\cite{MecklenbraukerLassoPath}.
For large $\mu$, the solution $\widehat{\mathbf x}={\mathbf 0}$ is trivial and $\mathbf{r} = 2\mathbf{A}^{H}\mathbf{y}$ in Eq.~\eqref{eq:ResidualVector}.
 Decreasing $\mu$, a first component of $\mathbf{x}$ becomes active when the corresponding component of $\mathbf{r}$ hits the boundary, $\mu = 2 \|\mathbf{A}^{H}\mathbf{y}\|_\infty$, ${\bf r}_i=\mu$. Inserting this solution into Eq.~\eqref{eq:ResidualVector} and solving for the second peak of ${\bf r}$ hitting the boundary $\mu$ indicates the value of $\mu$ for which a second component becomes active. This way we follow the LASSO path in Fig.\ \ref{fig:lassopath} towards less sparse solutions and lower $\mu$ as detailed in Ref.~ \onlinecite{MecklenbraukerLassoPath}.

Starting from Eq.\ (\ref{eq:ResidualVector}) with  $\Vec{x}_{\ell_1}(\mu^i)$ corresponding to regularization $\mu^i$ for the set of active indexes $\mathcal{M}_i$, the residual for the $n$th steering vector is now found.
\begin{eqnarray} 
{r}_{n}(\mu^i) &= &2 \Vec{a}_{n}^H \left( \Vec{y}- \Mat{A}\Vec{x}_{\ell_1}(\mu^i) \right) \nonumber \\
     &= & 2 \Vec{a}_{n}^H \left( \Vec{y}- \!\! \sum\limits_{m \in \mathcal{M}_i} \Vec{a}_m \Vec{x}_{{\ell_1},m}(\mu^i) \right)  \nonumber \\ 
     &\approx & 2 \Vec{a}_{n}^H \left( \Vec{y}- \sum\limits_{m \in \mathcal{M}_{i-1}} \Vec{a}_m \Vec{x}_{{\ell_1},m}(\mu^{i-1}) \right)  \label{eq:iterative_mu}\\  
      &\approx & 2 \Vec{a}_{n}^H  \Vec{y} \label{eq:fast_mu}
\end{eqnarray}
The two progressively stronger approximations, in Eqs.\ (\ref{eq:iterative_mu})--(\ref{eq:fast_mu}), above are valid if the steering vectors  corresponding to the final active set is sufficiently incoherent $|\Vec{a}_n^H \Vec{a}_m | \approx 0$. Eq.\ (\ref{eq:fast_mu}) actually corresponds to the conventional beamformer $\Vec{A}^H  \Vec{y}$ for a single snapshot. The above equation is  used for the selection of $\mu$ so it does not mean that the peaks in the conventional beamformer corresponds to the CS solution.

The procedure is given in Table \ref{algo:fast}, where ${\rm peak}({\bf r},k)$ is the $k$th peak of ${\bf r}$. We choose $F=0.9$.

The dual method has been used to estimate the solution path of the real-valued~\cite{TibshiraniLassoPath2011} and the complex-valued\cite{MecklenbraukerLassoPath} generalized  LASSO problems. The generalized LASSO  uses the $\ell_{1}$-norm to enforce structural or geometric constraints on the solution  by replacing the sparsity constraint  $\lVert \mathbf{x} \rVert_{1}$ with  $\lVert \mathbf{D}\mathbf{x} \rVert_{1}$ for a structured matrix $\mathbf{D}$. The generalized formulation  performs well in certain applications, e.g., recovery of continuous sources by promoting block sparsity~\cite{TibshiraniFusedLasso2005} and DOA tracking for moving sources by an adaptive update of a diagonal weighting matrix $\mathbf{D}$ which reflects the evolution of the source probability distribution~\cite{Mecklenbrauker:2013}.

\subsection{\label{sec:LassoPath}Regularization parameter selection via the LASSO path}

\begin{figure}[tb] 
\centering
\includegraphics[width=8.6cm]{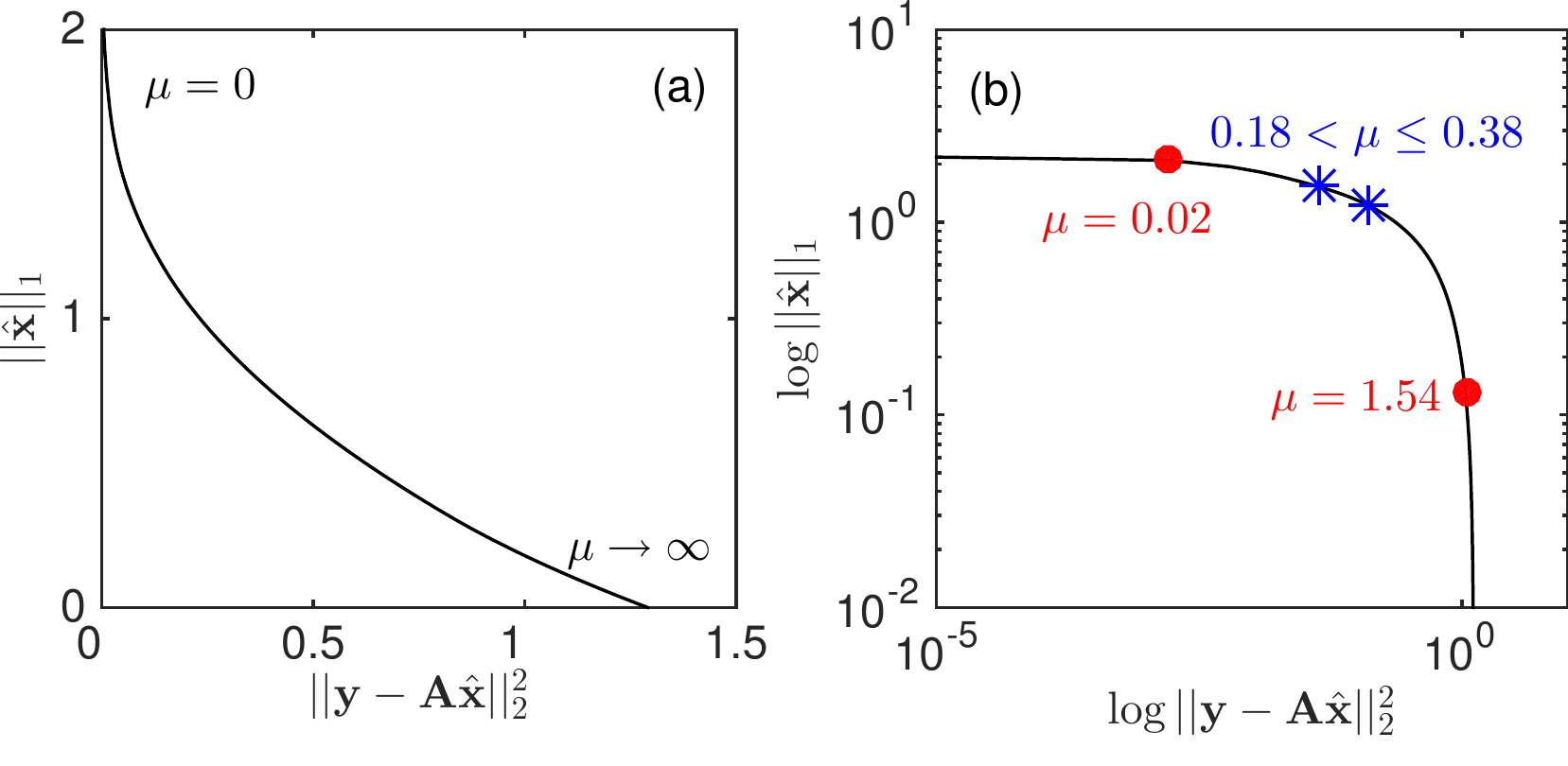}
\caption{(Color online) The data error $\lVert \mathbf{y} - \mathbf{A}\widehat{\mathbf{x}}\rVert_{2}^{2}$, describing the goodness of fit, versus the  $\ell_{1}$-norm in (a) linear scale and (b) log-log scale for the estimated solution $\widehat{\mathbf{x}}$ for different values of the regularization parameter $\mu$ in the LASSO problem Eq.~\eqref{eq:CS_solution_lasso} for sparse DOA estimation.}
\label{fig:Lcurve}
\end{figure}

\begin{figure}[tb]
\centering
\includegraphics[width=8.6cm]{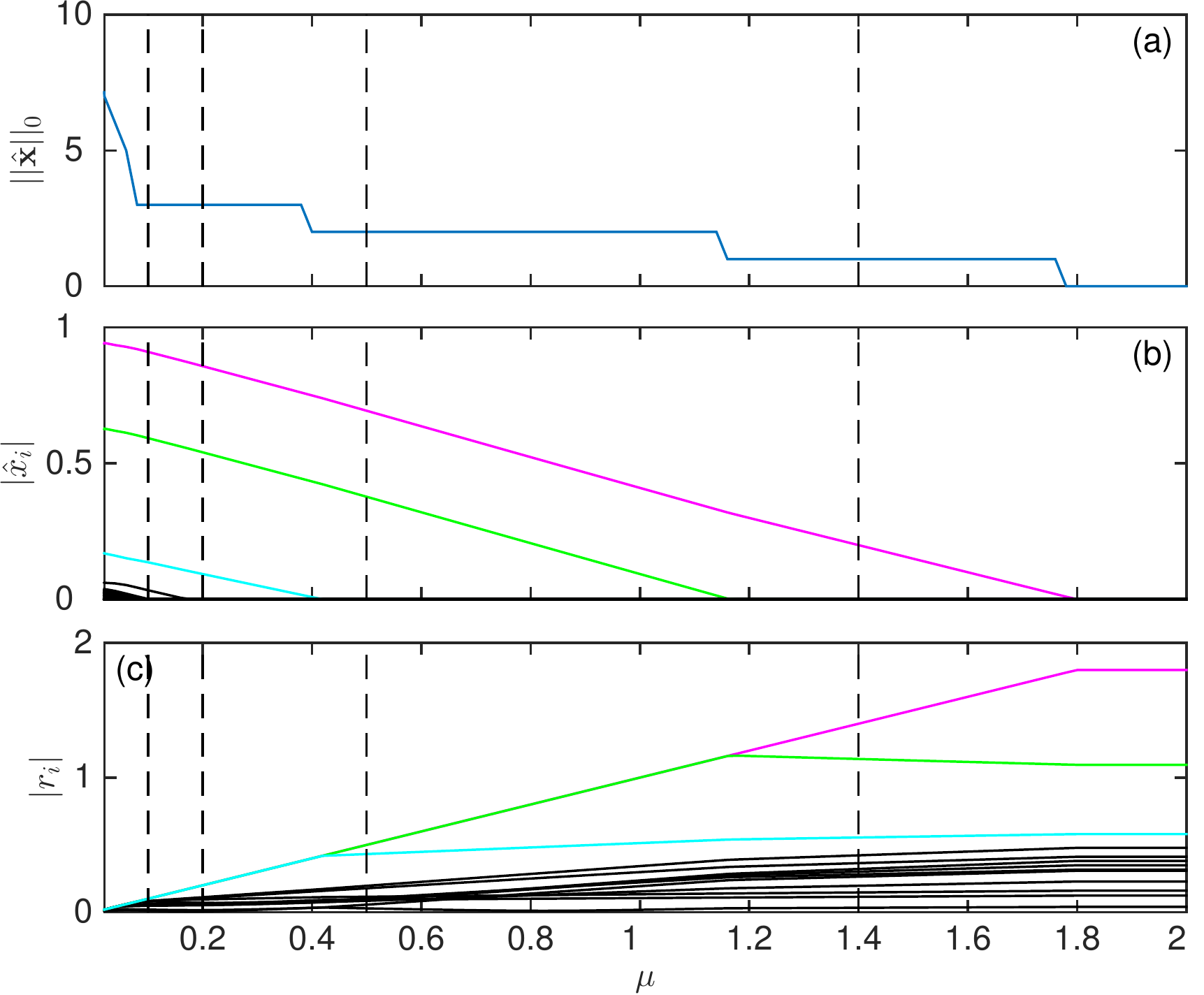}
\caption{(Color online) The LASSO path versus $\mu$ for three sources and SNR=20 dB. (a) Sparsity level of the estimate $\widehat{\mathbf{x}}$. (b) Paths for each component of the solution $\widehat{\mathbf{x}}$. (c) Paths for each component of the beamformed residual  $|\mathbf{r}| = 2|\mathbf{A}^{H}\left(\mathbf{y}- \mathbf{A}\widehat{\mathbf{x}} \right)|$. The vertical dashed lines indicates values of $\mu$ used in Figs.\ \ref{fig:LassoRegSolution} and \ref{fig:CohLassoPath}.}
\label{fig:LassoPath}
\end{figure}

\begin{figure}[tb]
\centering
\includegraphics[width=8.6cm]{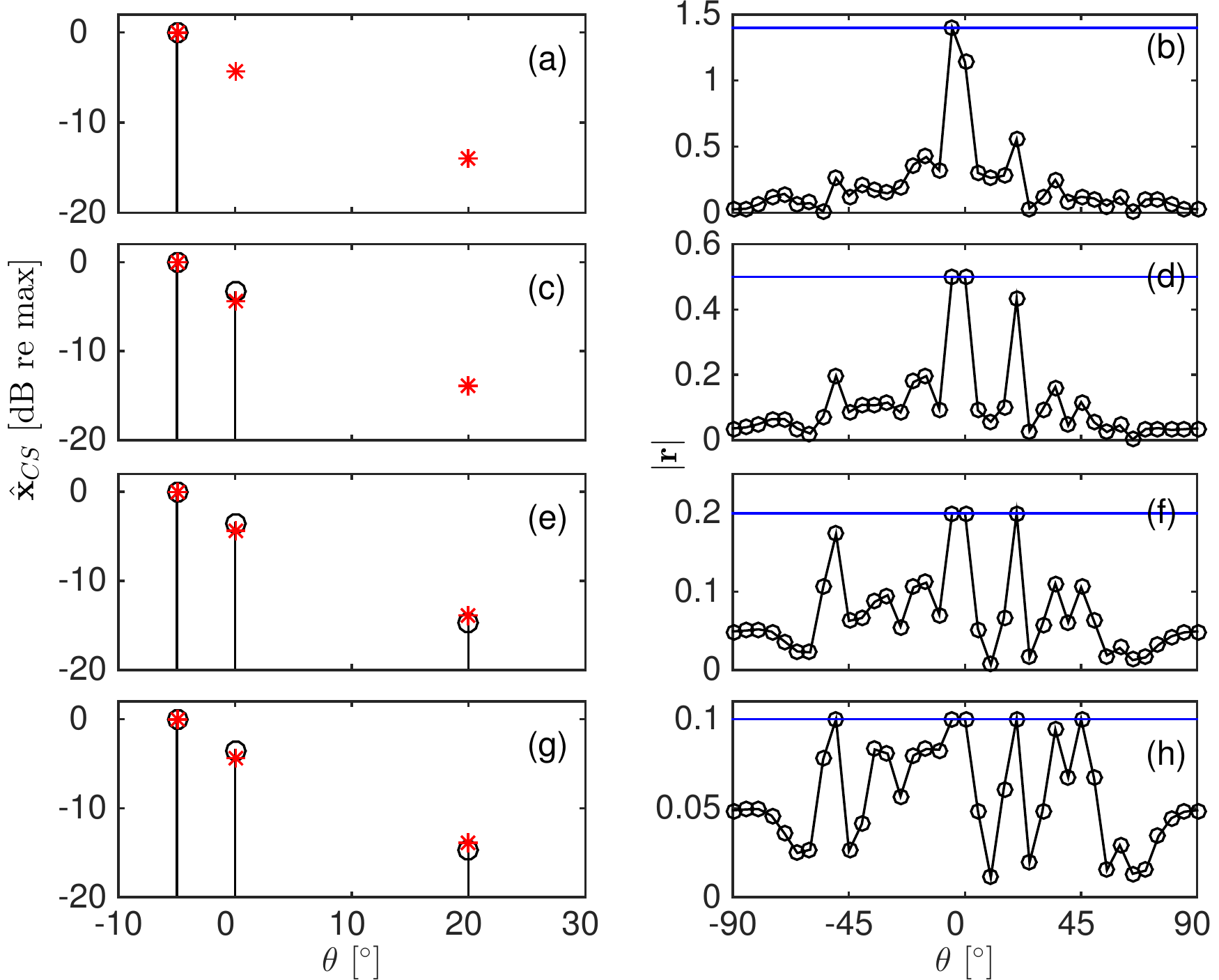}
\caption{(Color online) The unbiased estimate $\widehat{\mathbf{x}}_{\text{CS}}$ ($\circ$) for the true source $\mathbf{x}$ ($\star$) and the corresponding beamformed residual vector for  a coarse angular grid [$-90$:$5$:$90$]$^{\circ}$.   (a)--(b) $\mu = 1.4$, (c)--(d) $\mu = 0.5$, (e)--(f) $\mu = 0.2$, and  (g)--(h) $\mu = 0.1$   (corresponding to dashed lines in Fig.~\ref{fig:LassoPath}). The horizontal line in the residual plot (right) indicates the value of $\mu$.}
\label{fig:LassoRegSolution}
\end{figure}

\begin{figure}[tb]
\centering
\includegraphics[width=8.6cm]{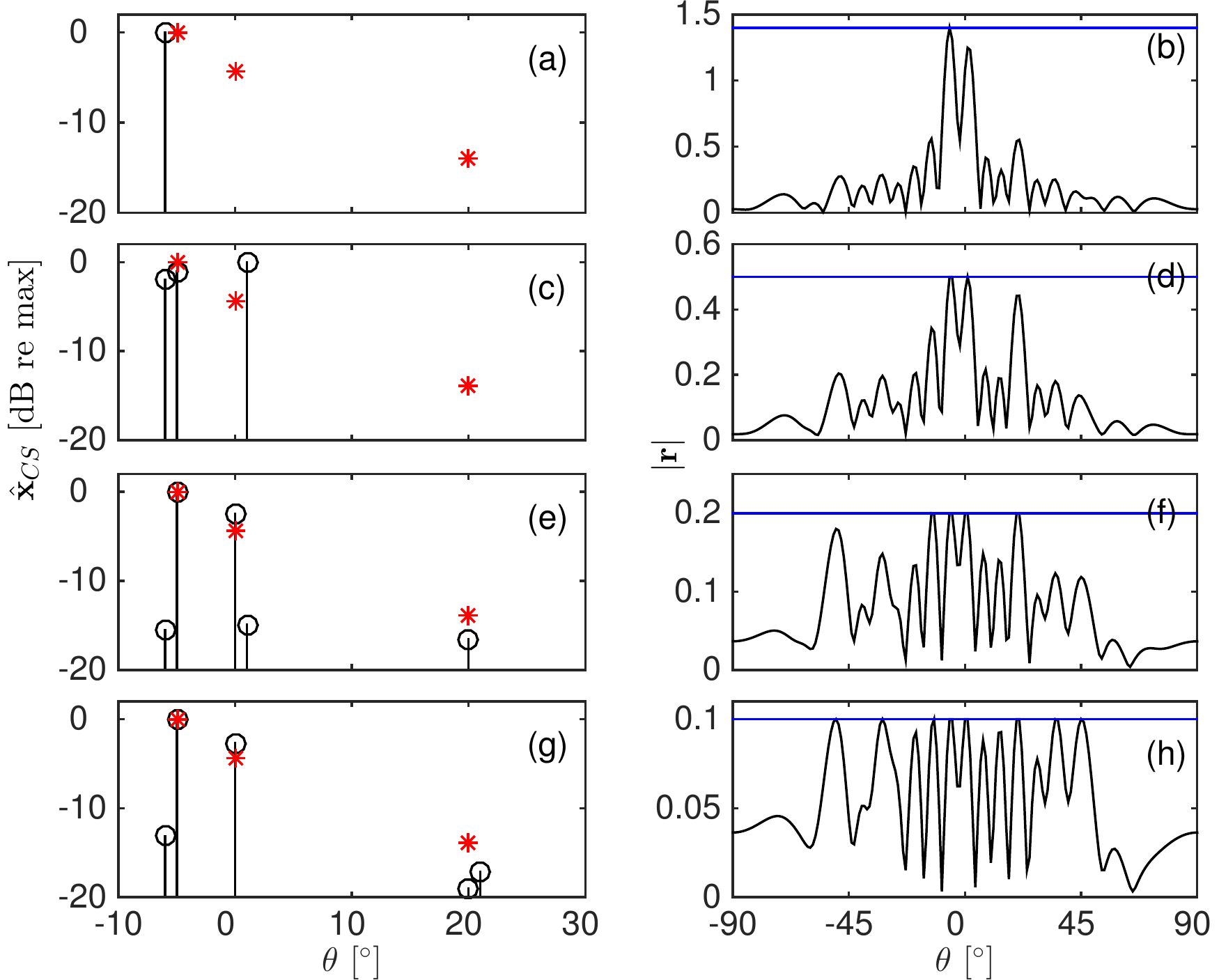}
\caption{(Color online) As in Fig.~\ref{fig:LassoRegSolution} but for the denser angular grid [$-90$:$1$:$90$]$^{\circ}$.  }
\label{fig:CohLassoPath}
\end{figure}

The LASSO performance in DOA estimation is evaluated by simulations starting with a large $\mu$ and subsequently decreasing its value. We consider an ULA with $M=20$  sensors and spacing $d = \lambda/2$.
Three sources are at DOAs $[-5,\, 0,\, 20]^{\circ}$ with corresponding magnitudes [1, 0.6, 0.2] (linear) or $[0, -4, -14]$ dB.
The sensing matrix $\mathbf{A}$ in (\ref{eq:SensingMatrix}) is defined on a coarse angular grid [$-90^{\circ}$:$5^{\circ}$:$90^{\circ}$] (Fig. \ref{fig:Lcurve}--\ref{fig:LassoRegSolution}) and a denser grid [$-90^{\circ}$:$1^{\circ}$:$90^{\circ}$] (Fig. \ref{fig:CohLassoPath}).
The noise variance in (\ref{eq:ArraySNR}) is chosen such that SNR=$20$ dB. 

The trade-off between regularization term $\lVert \widehat{\mathbf{x}}\rVert_{1}$  and the data fit $\lVert \mathbf{y} - \mathbf{A}\widehat{\mathbf{x}}\rVert_{2}^{2}$ in the LASSO estimate , Eq.~\eqref{eq:CS_solution_lasso}, for a range of values of $\mu$ is depicted in Fig.~\ref{fig:Lcurve}.  
The  relevant values of $\mu$ for the LASSO path are found between the two dots in Fig.~\ref{fig:Lcurve}(b), i.e. $1.54>\mu>0.02$. For these values of $\mu$, the importance shifts from favoring sparser solutions for large $\mu$ towards diminishing the model residual's $\ell_{2}$-norm for smaller $\mu$.
From inspecting  Fig.~\ref{fig:Lcurve}(b), it is difficult inferring the value of $\mu$ which results in  the desired sparsity level. The LASSO path offers a more insightful method to determine the range of good values of $\mu$ (contained within the asterisks in Fig.~\ref{fig:Lcurve}(b)) as explained below.

Figure~\ref{fig:LassoPath} shows (a) the sparsity level $\lVert\widehat{\mathbf{x}}\rVert_{0}$ of the LASSO solution, (b) the properties of the LASSO path and (c) the corresponding residual vector versus the regularization parameter $\mu$ Since the interest is on sparse solutions $\widehat{\mathbf{x}}$, it is  natural inspecting the LASSO path for decreasing values of $\mu$, i.e., interpreting Fig.~\ref{fig:LassoPath} from right to left.

For large values of $\mu$ (e.g., $\mu=2$) the problem in Eq.~\eqref{eq:CS_solution_lasso} is over-regularized, forcing  the trivial solution  $\widehat{\mathbf{x}} =\mathbf{0}$ (Fig.~\ref{fig:LassoPath}(b)), thus $\lVert\widehat{\mathbf{x}}\rVert_{0} =0$ (Fig.~\ref{fig:LassoPath}(a)). In this case, the slopes for all components $\lvert r_{i} \rvert$ are zero (Fig.~\ref{fig:LassoPath}(c)) since $\lvert r_{i} \rvert = \lvert 2 \mathbf{a}_{i}^{H}\mathbf{y}\rvert< \mu$ for all $i \in [0,\cdots,N]$ which is independent of $\mu$.

The first non-zero component of $\widehat{x}$ appears at $\mu \!\!=\!\! 2 \| \mathbf{A}^{H}\mathbf{y}\|_\infty \!\!=\!\! 1.76$ and remains active for $\mu \leq 1.76$ (Fig.~\ref{fig:LassoPath}(b)) increasing the sparsity level to $\lVert\widehat{\mathbf{x}}\rVert_{0} \!\!= \!\!1$ (Fig.~\ref{fig:LassoPath}(a)). The corresponding component $\lvert r_{i} \rvert \!\!= \!\!\lvert 2 \mathbf{a}_{i}^{H}\left(\mathbf{y}-\mathbf{a}_{i}\widehat{x}_{i}\right)\rvert$ (Fig.~\ref{fig:LassoPath}(c)) is equal to $\mu$ for $\mu\le 1.76$. The other components $r_{j}$ change slope at the singular point $\mu =1.76$, since now $\lvert r_{j} \rvert \!\!=\!\!\lvert 2 \mathbf{a}_{j}^{H}\left(\mathbf{y}-\mathbf{a}_{i}\widehat{x}_{i}\right)\rvert<\mu$ for all $j \in [0,\cdots,N]$, $j\neq i$. For $\mu\leq 1.14$, $\lVert\widehat{\mathbf{x}}\rVert_{0} = 2$ (Fig.~\ref{fig:LassoPath}(a)) as $\widehat{\mathbf{x}}$ acquires a second non-zero component (Fig.~\ref{fig:LassoPath}(b)) and the corresponding component $\lvert r_{i} \rvert$ becomes equal to $\mu$ (Fig.~\ref{fig:LassoPath}(c)). Similarly, the estimated solution has a third non-zero component for $\mu\leq 0.38$.

For $\mu\leq 0.18$, $\widehat{\mathbf{x}}$ has many non-zero components (Figs.~\ref{fig:LassoPath}(b),(c)) and its sparsity level increases abruptly (Fig.~\ref{fig:LassoPath}(a)). For such low values of $\mu$ the  importance shifts to the data fitting term ($\ell_{2}$-norm term)  in the regularized problem, Eq.~\eqref{eq:CS_solution_lasso}, and $\widehat{\mathbf{x}}$  includes many non-zero noisy components gradually reducing the data error.  

The specific values of $\mu$ at which an element of $\widehat{\mathbf{x}}$ becomes active are denoted as the singular points in the piecewise smooth LASSO path. 
At a singular point, some component of $\mathbf{r}$ hits the boundary $\mu$, i.e. $|r_n|=\mu$ for some index $n$.
Thus, the properties of the LASSO path indicate the selection of the regularization parameter $\mu$. For example, for a predefined sparsity level $K$ a good choice of $\mu$ is found by decreasing $\mu$ until the $K$th singular point at the LASSO path. 

Owing to the piecewise smooth nature of the LASSO path, there is a range of  $\mu$ which give  the same sparsity level (i.e., between two  singular points).  In principle, the lowest  $\mu$ in this range is desired as it gives the best data fit.  Though, any value of $\mu$ which achieves the desired sparsity suffices as once the active DOAs are recovered, the unbiased amplitudes are determined from Eq.~\eqref{eq:AmplitudeEstimation}.

Figure~\ref{fig:LassoRegSolution} shows the unbiased solution, Eq.~\eqref{eq:AmplitudeEstimation}, along with the corresponding beamformed residual for four sparsity levels  of $\mu$. Notice how the residuals decrease in value as $\mu$ is reduced.
For $\mu=0.1$,  Fig.\ \ref{fig:LassoRegSolution} shows that  five potential source locations exists as they have hit the boundary, so that components of $|{\bf r}|$ becomes equal to $\mu$. Solving
Eq.\ \eqref{eq:AmplitudeEstimation} shows that two sources are weak and are not shown. 

To increase precision in the LASSO reconstruction, a finer angular grid is required. However, angular grid refinement also causes higher coherence among steering vectors,  Eq.~\eqref{eq:SteeringVector}, and the problem in Eq.~\eqref{eq:NoisyData} becomes increasingly underdetermined. Then,  when solving the LASSO minimization Eq.~\eqref{eq:CS_solution_lasso} might not exhibit the desired sparsity. 
Due to basis coherence and as $\mu$ decreases, components in the estimate $\widehat{\mathbf{x}}$ can be either activated (become non-zero)  or annihilated.
Similarly   the residual components can hit  or leave the boundary\cite{TibshiraniLassoPath2011} (where components of $ | {\bf r} |$ is equal to $\mu$, see Eq.\ (\ref{eq:BoundedResidual})).

In Fig.~\ref{fig:CohLassoPath}  a denser angular grid with spacing $1^{\circ}$ and setup as  Fig.~\ref{fig:LassoRegSolution} is used.
For $\mu=1.4$, there is just one active component (Fig.~\ref{fig:CohLassoPath}(a)) at $-6^{\circ}$ which is $1^{\circ}$ away from the strongest DOA. 
This offset is mainly due to basis coherence, the correct location is not yet recovered. As $\mu$ is decreased the correct bin is eventually obtained (Fig.~\ref{fig:CohLassoPath}(e)). 
Thus when searching for a $K$ sparse solution, it is often advantageous to search initially for more than $K$ peaks, 
and then limit the final solution to the $K$ most powerful elements.

The residual $\bf r$ is systematically decreased as $\mu $ is reduced. All the active component can be seen as where  components of $|{\bf r}|$ becomes equal to $\mu$ in the right panel of  Fig.~\ref{fig:CohLassoPath}.

 \section{DOA estimation error evaluation} 
\label{sec:perf}

If the source DOAs are well separated with not too different magnitude, the DOA estimation for multiple sources using CBF and CS turns out to behave similarly. 
They differ, however, in their behavior whenever two sources are closely spaced.  The same applies for MVDR under the additional assumptions of incoherent arrivals and sufficient number of snapshots, $L\geq M$.
The details are of course scenario dependent.

For the purpose of a quantitative performance evaluation with synthetic data, the estimated, $\widehat\theta_k$, and the true, $\theta_k^{\rm true}$, DOAs are paired with each other such that the root mean squared DOA error is minimized in each single realization.
After this pairing, the ensemble root-mean-squared error  is computed,
\begin{equation}
\mathrm{RMSE}=\sqrt{{\rm E}\left[\frac{1}{K}\sum_{k=1}^K(\widehat\theta_k-\theta_k^{\rm true})^2\right]}~.
\end{equation}
%

The data is generated to have a fixed SNR Eq.\ \eqref{eq:ArraySNR}. The source phases of each $\bf x$ component  is uniformly distributed on $[0, 2\pi )$ in order to generate a sample covariance matrix from which  MVDR/MUSIC can resolve incoherent sources.

CBF  suffers from low-resolution and the effect of sidelobes for both single and multiple data snapshots, thus the simple peak search used here is too simple. These problems are reduced in MVDR and MUSIC for multiple snapshots and they do not arise with CS.

The optimal performance for $K$ sources is found by searching over all combinations of steering vectors for the maximum likelihood solution, Eq. \eqref{eq:Likelihood}, i.e., the best fitting source vector using Eq.~(\ref{eq:AmplitudeEstimation}). This is a NP-hard combinatorial problem, that for $N$ look directions requires evaluation of $N!/K!(N-K)!$ solutions. For $N\!\!=\!\!361$ and either $K\!\!= \!\!2$ or  $K\!\!=\!\! 3$, this gives 77,000 or 7,700,000 combinations to be evaluated. This makes the exhaustive search approach impractical for larger $K$.

\begin{figure}[tb] 
   \includegraphics[width=0.45\textwidth]{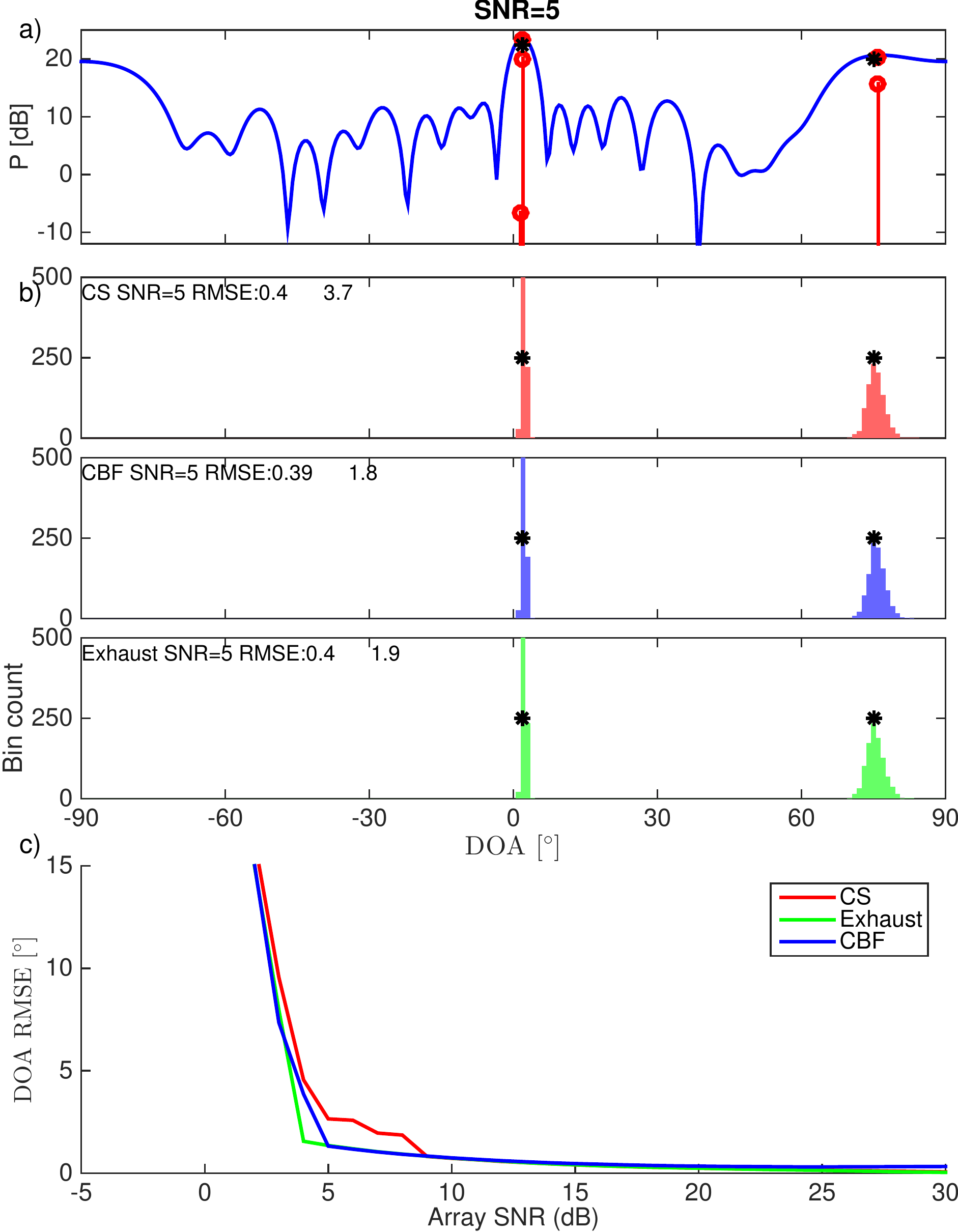}
       \caption{ (Color online) Single snapshot example for 2 sources at DOAs $[2,\,75]^\circ$ and magnitudes $[22 ,\,20]$ dB. At $\mathrm{SNR}=5\,\mathrm{dB}$ a) spectra for CBF, CS (o) and unbiased CS (o, higher levels), and b) CS, CBF and exhaustive-search histogram based on 1000 Monte Carlo simulations,  and c) CS, CBF and exhaustive-search performance versus SNR. The true source positions (*) are indicated in a) and b).
    \label{fig:SMV2}}
\end{figure}

\begin{figure}[tb] 
   \includegraphics[width=0.45\textwidth]{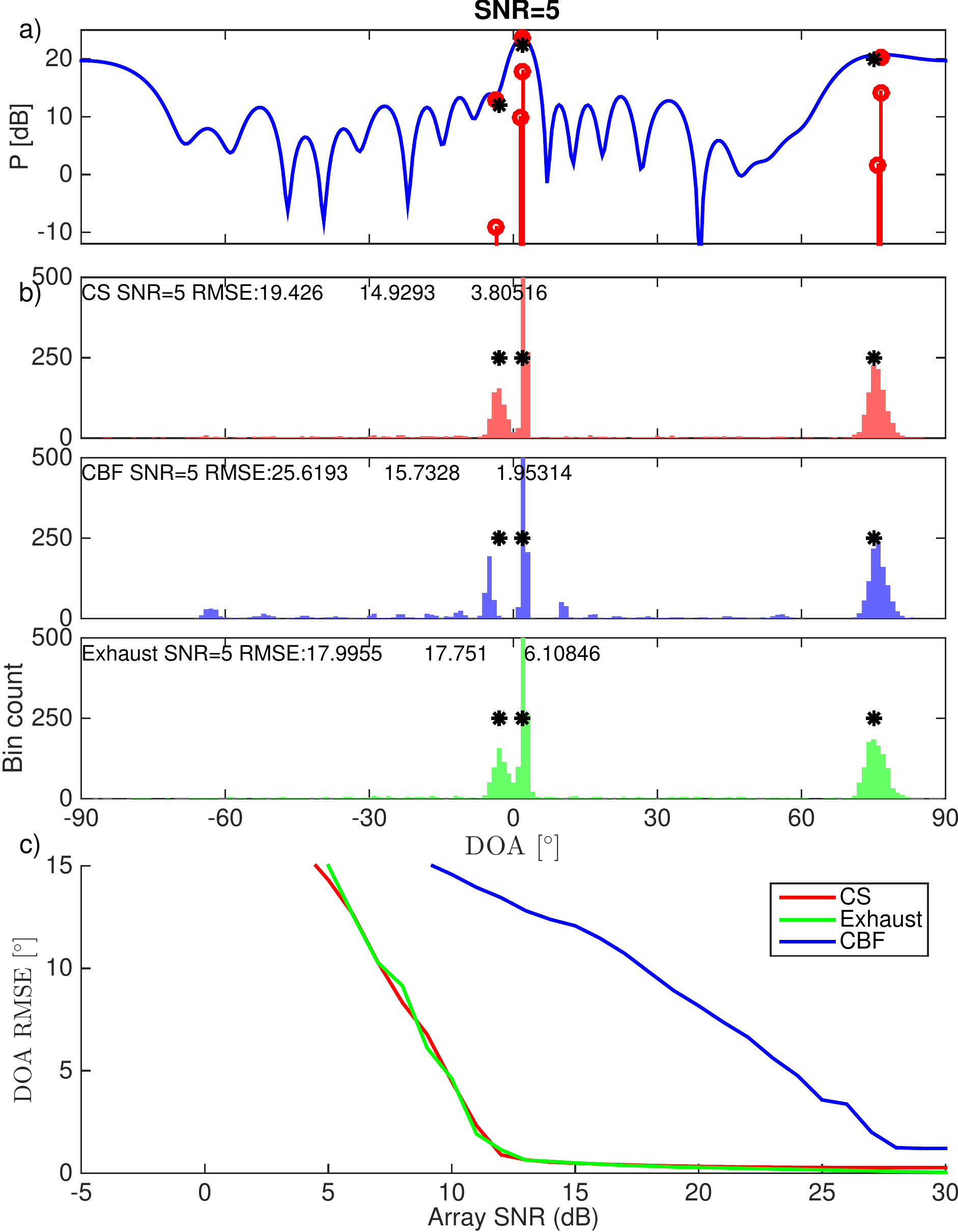}
    \caption{ (Color online) As Fig.\ \ref{fig:SMV2} but for 3 sources at DOAs $[-3,\,2,\,75]^\circ$ and magnitudes $[12,\,22,\,20]$ dB.
       \label{fig:SMV}}
\end{figure}

In the following simulation, we consider an array with $M=20$ elements with spacing $d=\lambda/2$. The DOAs  are assumed to be on a fine angular grid [$-90^{\circ}$:$0.5^{\circ}$:$90^{\circ}$], i.e.\ $\mathbf{A} \in \mathbb{C}^{20\times 361}$. 
The regularization parameter $\mu$ is chosen to correspond to the $K+2$ largest peak of the residual in Eq.~\eqref{eq:ResidualVector} using the procedure in Table \ref{algo:fast} and retaining only the K largest source powers.
We require the peaks of the CS to be at least 4 bins apart. Thus the exhaustive and the CS do not solve the identical problem, as the CS solves a smaller problem.
Note that panel c in Figs.  \ref{fig:SMV2}--\ref{fig:MMV2}  shows the simulation results versus 
array SNR defined in Eq.~\eqref{eq:ArraySNR}.


\subsection{Single Snapshot}
In the first scenario, we consider a single snapshot case with additive noise 
  with $K=2$ well-separated DOAs at $[2,\,75]^\circ$ with magnitudes $[22,\,20]$ dB, see Fig. \ref{fig:SMV2}. 
In the second scenario, a third weak source is included very close to the first source:
Thus, $K=3$ and the source DOAs are $[-3,\, 2,\, 75]^\circ$ with magnitudes $[12,\, 22,\, 20]$ dB, see Fig. \ref{fig:SMV}. 
The synthetic data is generated according to Eq.~\eqref{eq:NoisyData}.

For the first scenario, the CS diagrams in  Fig. \ref{fig:SMV2}a show DOA estimation with small  variance but 
indicate a bias towards endfire, as  for the  true DOA $75^{\circ}$ the CS estimate is $76^{\circ}$.
Towards endfire the main beam becomes broader and absorbs more noise power,
The CBF spectra Fig. \ref{fig:SMV2}a  are characterized by a high sidelobe level but for  the two well-separated similar-magnitude sources this is a minor problem. 

Using Monte Carlo simulations, we repeat the CS inversion for 1000 realizations of the noise in Fig. \ref{fig:SMV2}b. The RMSE  increases towards the endfire directions. This is to be expected as the main beam becomes wider and this results in a lower DOA resolution \cite{XenakiCS:2014}. Since the sources are well-separated in this scenario, CS, CBF, and exhaustive search perform similarly with respect to RMSE.

Repeating the Monte Carlo simulations at several SNRs  gives the RMSE performance of CS and CBF in Fig. \ref{fig:SMV2}c.
Their performance is about the same since the DOAs are well-separated.

In the second scenario, the CBF cannot resolve the two closely spaced sources with DOAs $[-3,\, 2]^\circ$. 
They are less than a beamwidth apart as indicated in Fig. \ref{fig:SMV}a.
Sidelobes cause a few DOA estimation errors at $-65^\circ$ in the CBF histogram, Fig. \ref{fig:SMV}b. 
Since CS obtains high-resolution even for a single snapshot, it performs much better than CBF, Fig. \ref{fig:SMV}c. 

 \subsection{Multiple Snapshot}
  \label{se:mmv}

In the multiple-snapshot scenario,  MVDR and MUSIC use the data sample covariance matrix Eq.~\eqref{eq:SCM} whereas CBF and CS works directly on the observations $\mathbf X$ Eq.~\eqref{eq:CS_lasso_MultiSnap}.
The sample covariance matrix  is formed by averaging $L$ synthetic data snapshots.
The source magnitude is considered invariant across  snapshots. 
The source phase is for each snapshot sampled from a uniform distribution on $[0,\,2\pi)$.

Due to the weak performance of MVDR in scenarios with coherent arrivals\cite{VanTreesBook}, we assume incoherent arrivals in the simulations although not needed for CS.
For CS we use Eq.~\eqref{eq:CS_lasso_MultiSnap} with a similar choice of regularization parameter $\mu$ as for the single snapshot case. 

\begin{figure}[htb] 
   \includegraphics[width=0.45\textwidth]{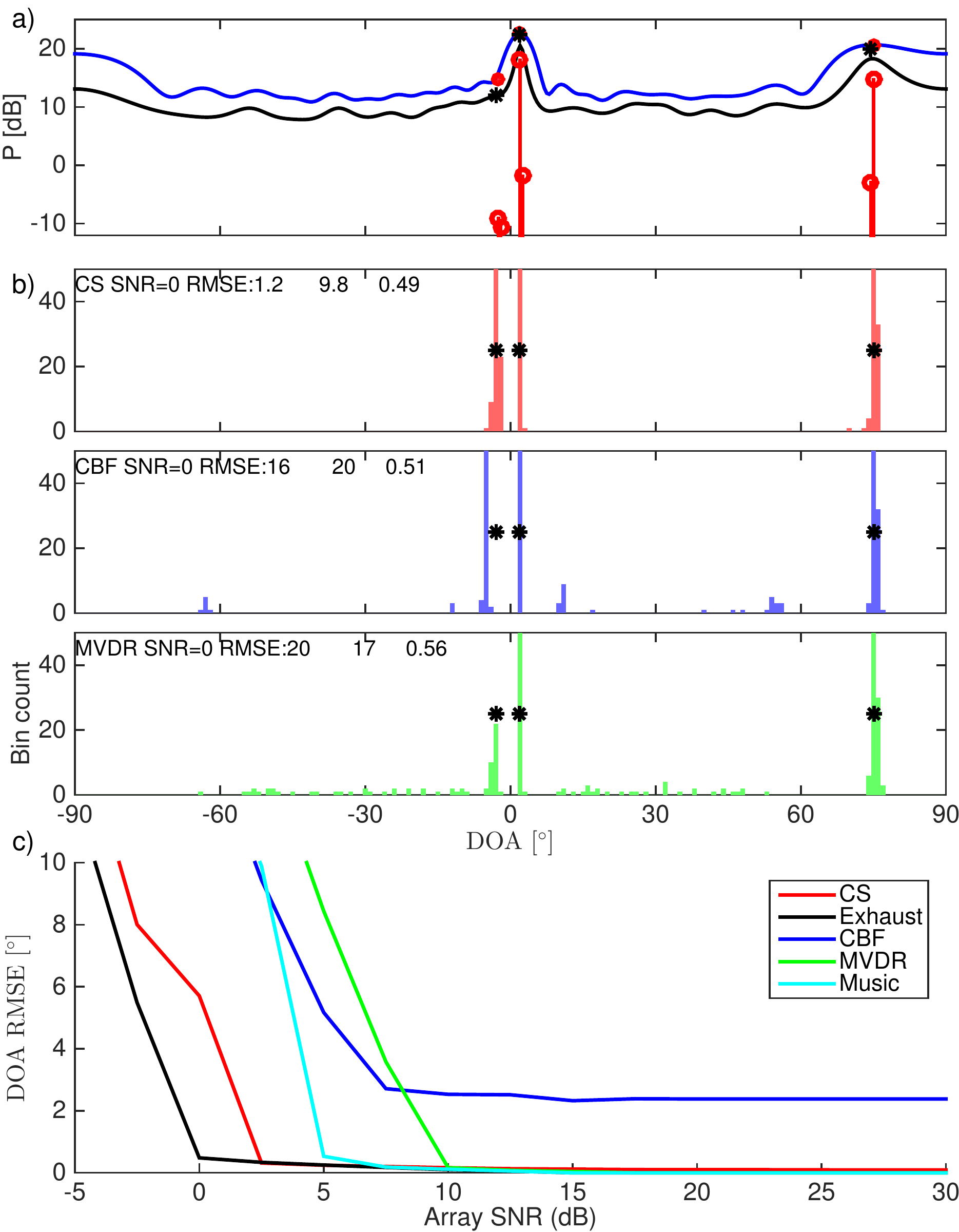}
    \caption{ (Color online) Multiple $L=50$ snapshot example for 3 sources at DOAs $[-3,\, 2,\, 75]^\circ$ with magnitudes $[12,\, 22,\, 20]$ dB. At $\mathrm{SNR}=0\,\mathrm{dB}$ a) spectra for CBF, MVDR, and CS (o) and unbiased CS (o, higher levels), and b) CS, CBF and MVDR histogram based on 100 Monte Carlo simulations,  and c)  CS,  exhaustive-search, CBF, MVDR, and MUSIC performance versus SNR. The true source positions (*) are indicated in a) and b).    \label{fig:MMV}}
\end{figure}

Using the same setup as in Fig.\ \ref{fig:SMV}, but estimating the source DOAs based on $L=50$ snapshots gives the results in Fig. \ref{fig:MMV}.
At $\mathrm{SNR}=0\,\mathrm{dB}$ the diagrams in Fig. \ref{fig:MMV}a show that CS localizes the sources well, in contrast to 
the CBF and MVDR that is also indicated in the  histograms in Fig. \ref{fig:MMV}b.
The RMSE in Fig. \ref{fig:MMV}c, shows that CBF does not give the required resolution even for high SNR.
MVDR performs well for $\mathrm{SNR}>10\,\mathrm{dB}$, whereas CS performs well for SNRs down to  $2.5\,\mathrm{dB}$.

\begin{figure}[htb] 
    \includegraphics[width=0.45\textwidth]{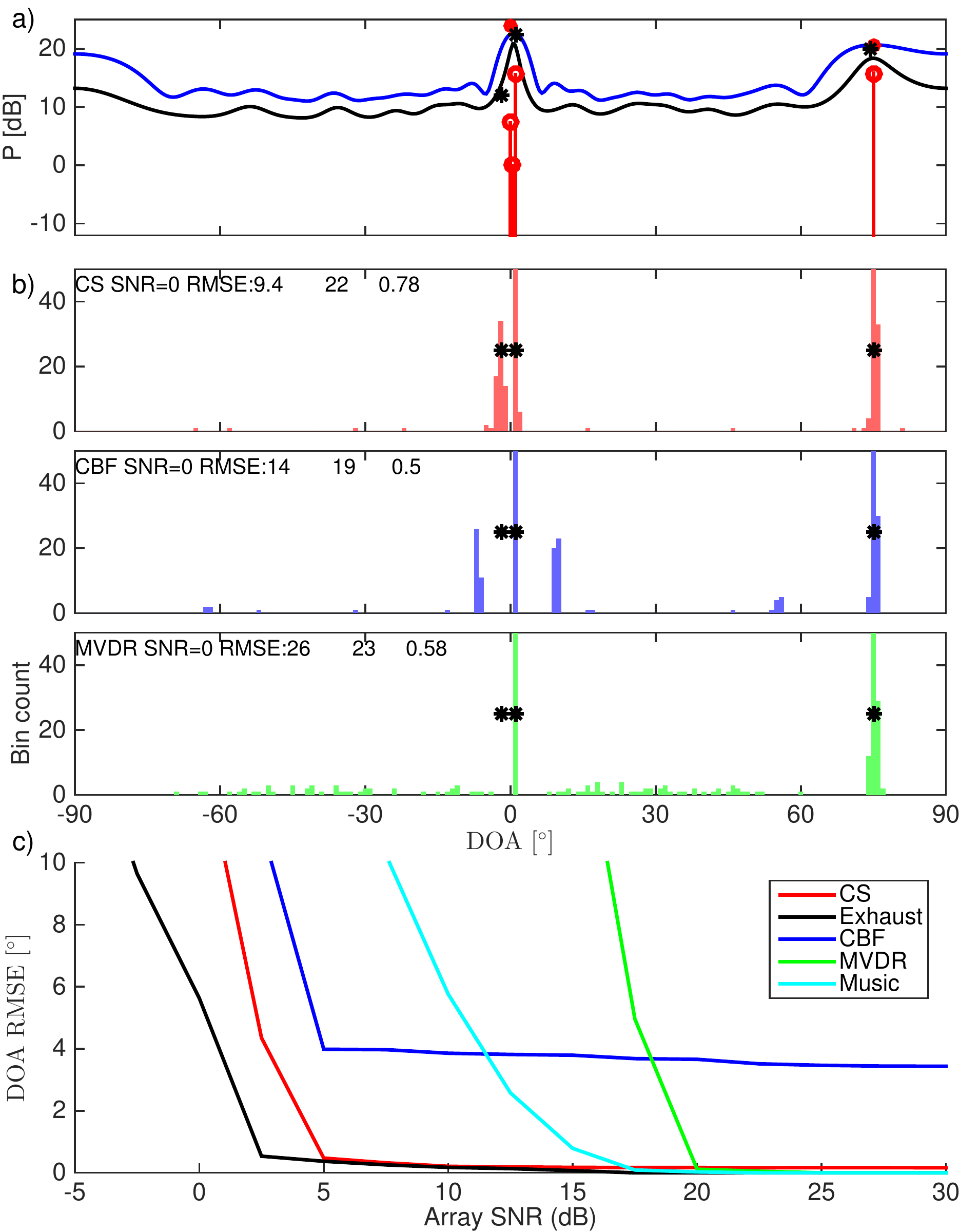}
    \caption{ (Color online) As Fig.\ \ref{fig:MMV} but with closer spaced sources [$-$2, 1, 75]$^\circ$. 
    \label{fig:MMV2}}
\end{figure}

In a third scenario, the weak broadside sources are moved closer with DOAs defined as $[-2,\, 1,\, 75]^\circ$. 
Fig. \ref{fig:MMV2} gives about the same DOA estimates for CBF, as it is already at its maximum performance even for high array SNR, confirming its low resolution. 
MVDR fails for $\mathrm{SNR}<20\,\mathrm{dB}$, which is $10\,\mathrm{dB}$ higher than the corresponding value in Fig. \ref{fig:MMV}c (MUSIC fails also at a level $10\,\mathrm{dB}$ higher). Contrarily, CS fails only for $\mathrm{SNR}<5\,\mathrm{dB}$ which is $2.5\,\mathrm{dB}$ higher  (Figs.~\ref{fig:MMV}c and~\ref{fig:MMV2}c).
Note how MVDR completely misses the weak source at $-2^\circ$ in Figs. \ref{fig:MMV2}c, but CS localize it with a larger spread.
%
Thus, as the weak source moves closer to the strong source, CS  degrades slower than MVDR in terms of RMSE. 
This is a good indication of its high-resolution capabilities.

\begin{figure}[htb] 
     \includegraphics[width=0.45\textwidth]{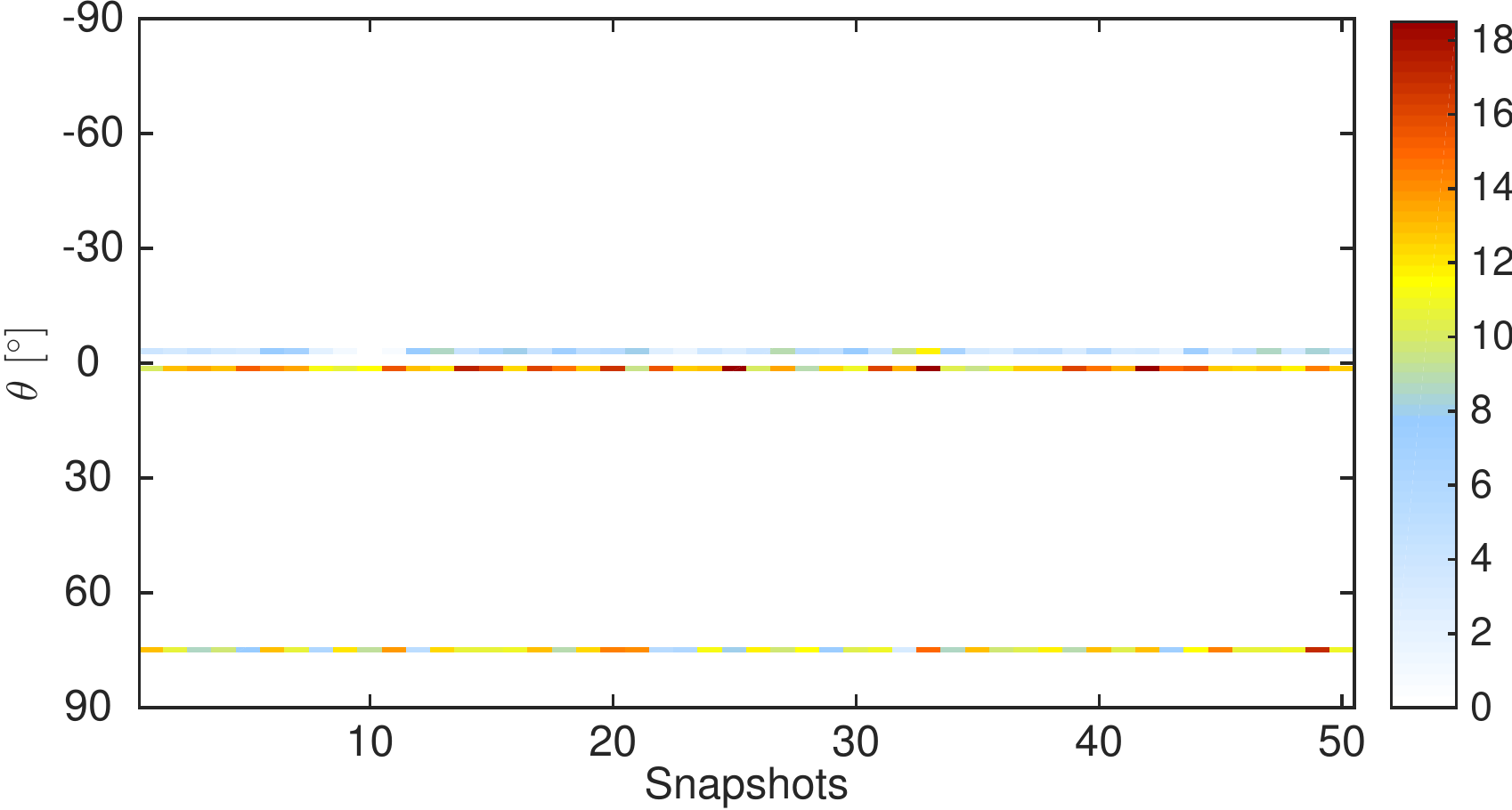}
    \caption{ (Color online)
 Power (linear) for the multiple snapshot case across azimuths and snapshots for one noise realization at 
  $\mathrm{SNR}=0\,\mathrm{dB}$ for the scenario with DOAs at $[-3,\, 2,\, 75]^\circ$.
    \label{fig:L1power}}
\end{figure}

Figure \ref{fig:L1power} shows the estimated power at the one realization in Fig.\ \ref{fig:MMV}a of $L=50$ snapshots inverted simultaneously. 
We emphasize the scale of the problem. Equation  \eqref{eq:MultiSnapNoisyData} has $20\cdot 50 = 1000$ equations to 
determine $361\cdot 50 =18050$ complex-valued variables at 361 azimuths and 50 snapshots observed on 20 sensors. 
The sparsity constraint is crucial here.

The CS (and especially the exhaustive-search) requires several orders of magnitude more CPU-time than the beamforming methods.

Many other simulations could be performed,  for example colored noise, no assumptions on number of sources, and random source locations. From initial exploration of these it is our impression that CS will perform well, though more simulations are required.


\begin{figure}[htb] 
\centering
\includegraphics[width=0.45\textwidth]{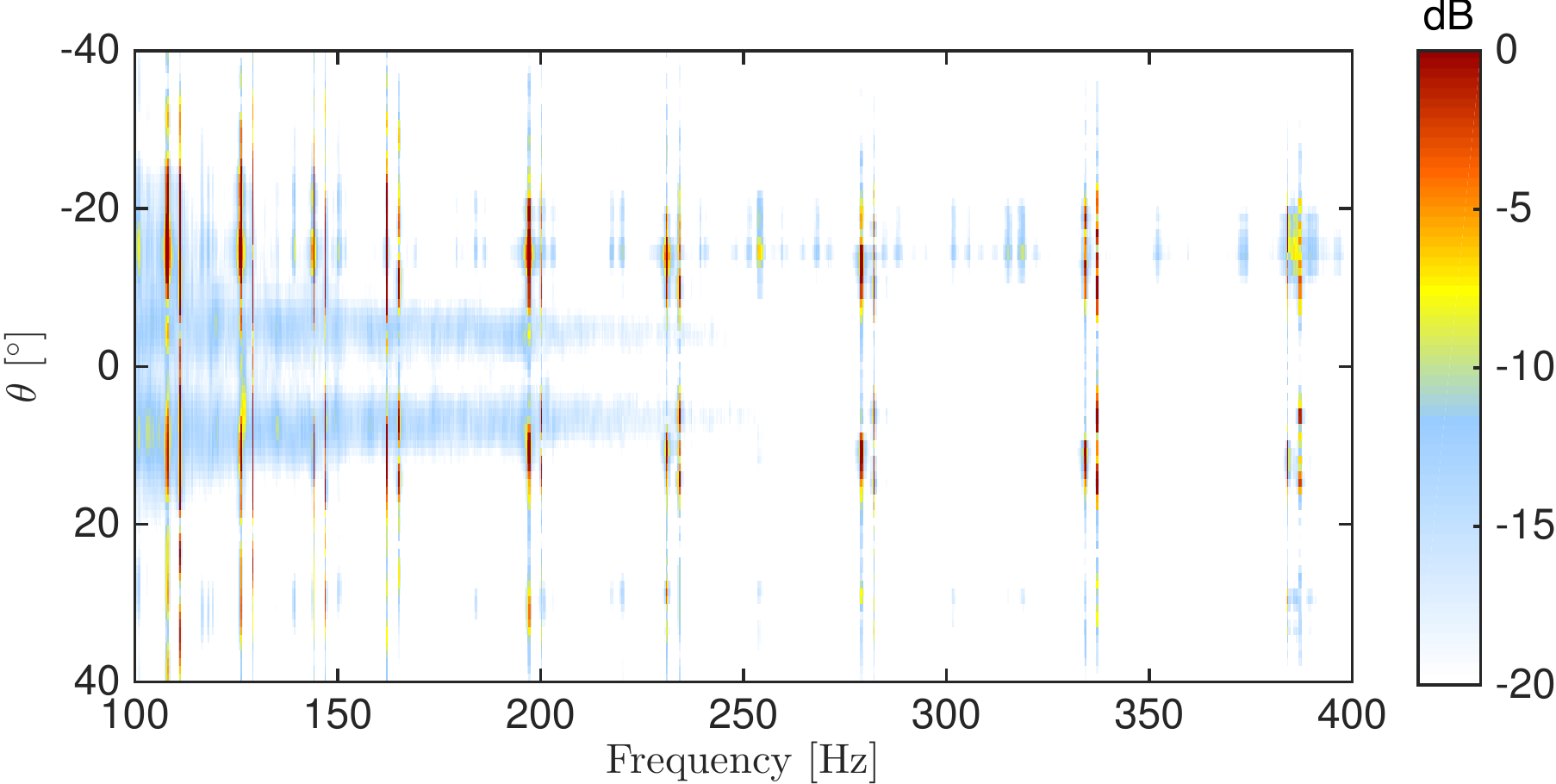}
\caption{(Color online) Spatial CBF spectrum across frequency at the source's closest point of approach to the array.}
\label{fig:beamfreq}
\end{figure}
\begin{figure*}[htb] 
\centering
\includegraphics[width=1\textwidth]{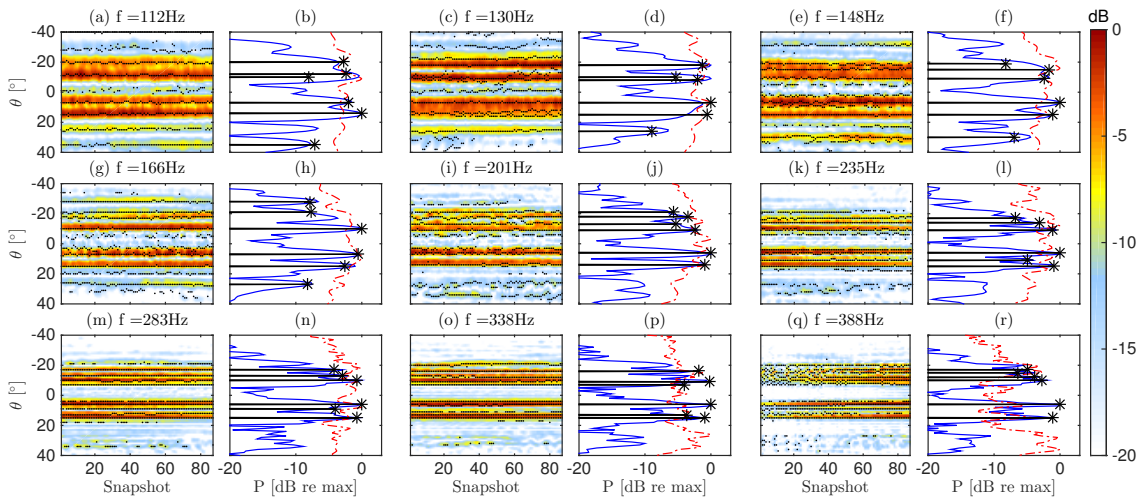}
\caption{(Color online) Single (contour plots) and multiple (line plots) snapshot reconstruction at the transmitted frequencies with CS ($\star$), CBF (background color, solid) and MVDR (dashed). For the single snapshot we have assumed $K=10$ sources while for the multiple snapshot $K=6$.}
\label{fig:deep}
\end{figure*}

\section{Experimental results}
\label{sec:swellex96}

The high-resolution performance of CS both in single- and multiple-snapshot cases is validated with experimental data in a complex multi-path shallow-water environment and it is compared with conventional methods, namely CBF and MVDR.

The data set is from the shallow water evaluation cell experiment 1996 (SWellEx-96) Event S5\cite{booth1996,dspain1999} collected on a 64-element vertical linear array. The array has uniform intersensor spacing 1.875~m and was deployed at  waterdepth 16.5~m spanning  94.125--212.25~m. During the Event S5, from 23:15--00:30 on 10-11 May 1996 west of Point Loma, CA, two sources, a shallow and a deep, were towed simultaneously from 9~km southwest to 3~km northeast of the array at a speed of 5 knots (2.5~m/s). Each source was transmitting a unique set of tones.

Here, we are interested in the deep source towed at 54~m depth while at the vicinity of the closest point of approach (CPA) which was 900~m from the array and occurred around 00:15, 60~min into the event. The deep-towed source signal submitted a set of 9 frequencies [112, 130, 148, 166, 201, 235, 283, 338, 388]~Hz at approximately 158~dB \textit{re} $1\mu\mathrm{Pa}$. The processed recording has duration of 1.5~min (covering 0.5~min before and 1~min after the CPA) sampled at 1500~Hz. It was  split into 87 snapshots of $2^{12}$ samples (2.7 s) duration, i.e., with  63\% overlap.

Figure~\ref{fig:beamfreq} shows the multiple-snapshot CBF spatial spectrum, Eq.~\eqref{eq:beam}, over the 50-400 Hz frequency range. Arrivals are detected not only at the transmitted tonal frequencies of the deep towed source but also at several other frequencies corresponding to the shallow-towed source tonal frequencies, weaker deep source frequencies, and the acoustic signature of the tow-ship.

Single-snapshot processing with CBF and CS at the deep source tonal set, contour plots in Fig.~\ref{fig:deep}, indicates the presence of several multipath arrivals which are adequately stationary along the snapshots at the CPA. Due to the significant sound speed variation it is not straightforward to associate the reconstructed DOAs with specific reflections. The CBF map comprises 6 significant peaks but suffers from low resolution and artifacts due to sidelobes and noise. To choose the regularization parameter in the LASSO formulation for CS reconstruction, we solve iteratively Eq.~\eqref{eq:CS_solution_lasso} as described in Table ~\ref{algo:fast} with initial value $\mu = 2\| \mathbf{A}^{H}\mathbf{y}\|_\infty$, until the obtained estimate has a sparsity level of 10. The CS reconstruction results in improved resolution due to the sparsity constraint and significant reduction of artifacts in the map.

Combining the data from all the snapshots and processing with CBF, MVDR, and CS, as in Sec \ref{se:mmv}, reveals that MVDR fails to detect the coherent multipath arrivals; see line plots in Fig.~\ref{fig:deep}. Again the peaks of CBF and CS are consistent but CS offers  improved resolution. 

We have here used higher sparsity for the single-snapshot processing to allow for identifying non-stationary paths. The non-stationary path can be seen in several of the contour plots, most prominently at 112, 130 and 201Hz.
When performing multiple-snapshot processing where the solution is constrained to remain active at one azimuth (but with varying power), the stationary paths are most likely to contribute to the CS solution.

\section{Conclusion}

The estimation of multiple directions-of-arrival (DOA) is formulated as a sparse source reconstruction problem. 
This is efficiently solvable using compressive sensing (CS) as a least squares  problem regularized with a sparsity promoting constraint. 
The resulting solution is the maximum a posteriori (MAP) estimate for both the single and multiple-snapshot formulations.  
The regularization parameter balances the data fit and the solution's sparsity. 
It is selected so that the solution is sufficiently sparse providing high-resolution DOA estimates. 
A procedure to find an adequate choice for the regularization parameter is described whereby the DOAs are obtained.

CS provides high-resolution acoustic imaging both with single and multiple snapshot.  
The performance evaluation shows that for single snapshot data, CS gives higher resolution than CBF.
For multiple snapshots, CS provides higher resolution than MVDR/MUSIC and the relative performance improves as the source DOAs move closer together.

The real data  example indicates that CS is capable of resolving multiple coherent wave arrivals, e.g. stemming from multipath propagation.
 
 \section*{Acknowledgment}
This work was supported by 
the Office of Naval Research, Grant Nos. N00014-11-1-0439 and N00014-13-1-0632 (MURI), as well as   FTW Austria's ``Compressed channel state information feedback for time-variant MIMO channels''.

\bibliographystyle{unsrt}
\bibliography{CSlasso}

\end{document}